\documentclass[a4paper,dvipsnames]{article}
\usepackage[english]{babel}
\usepackage{amsmath}
\usepackage{graphicx}
\usepackage{parskip} 
\usepackage{subcaption} 
\usepackage{multirow}	
\usepackage{mathtools} 
\usepackage{float} 
\usepackage{siunitx}
\usepackage{amssymb} 
\usepackage{comment}
\usepackage{xcolor}
\usepackage{todonotes}

\usepackage[noend]{algpseudocode}
\usepackage{algorithm,algorithmicx}

\usepackage{Bert_ii_variables_copy}

\usepackage{pgfplotstable}
\usepackage{pgfplots}
\usepackage{authblk}

\usepackage{tikz}
\usetikzlibrary{external,calc,arrows.meta}
\tikzexternaldisable

\mathtoolsset{showonlyrefs}

\definecolor{azure}{rgb}{0.0, 0.5, 1.0}
\usepackage{tcolorbox}
	\tcbuselibrary{most}

\newtcolorbox[auto counter, number within=section]{notebox}[1][]{%
	outer arc=0pt,
	arc=0pt,
	colframe=azure,
	colback=azure!20,
	attach boxed title to top left={yshift=-3mm, xshift=3mm},
	enhanced,
	boxed title style={
		colback=azure,
		top=3pt,
		bottom=3pt,
	},
  	fonttitle=\sffamily,
	top=12pt,
  	title=Note~\thetcbcounter,
}

\newcommand{\VAscaleparam}{\epsilon}

\renewcommand{\VArate}{{\color{red}\xi^{-1}}}

\renewcommand{\VAstandardexp}{\varepsilon}

\newcommand{\VAewstdv}{s}
\newcommand{\VAewmean}{m}

\title{Kinetic-diffusion asymptotic-preserving Monte Carlo algorithm for Boltzmann-BGK in the diffusive scaling}

\author[1,$\star$]{Bert Mortier}
\author[2]{Martine Baelmans}
\author[1]{Giovanni Samaey}
\affil[1]{KU Leuven, Department of Computer Science}
\affil[2]{KU Leuven, Department of Mechanical Engineering}
\affil[$\star$]{Corresponding author: Bert Mortier, bert.mortier@cs.kuleuven.be}

\date{}

\newcommand{\VAtotaltime}{\bar{t}}
\newcommand{\VAeventtime}{\tau}
\newcommand{\VAtimevar}{t}
\newcommand{\VAtimestep}{\Delta \VAtimevar}
\newcommand{\VAeventtimedif}{\Delta \VAeventtime}

\newcommand{\VAtimestepno}{n}
\newcommand{\VAnoftimesteps}{N}

\newcommand{\VAeventnointimestep}{{\VAeventno}}

\newcommand{\VAnofeventsintimestepnon}{{\VAnofevents}}

\newcommand{\VAeventnotilltimestep}[1]{{\kappa^{#1}}}

\newcommand{\VAhermite}[1]{H_{#1}}

\newcommand{\VAexpec}[1]{\mathbb{E}\!\left[#1\right]}
\newcommand{\VAexpecover}[2]{\mathbb{E}_{#1}\!\left[#2\right]}
\newcommand{\VAvar}[1]{\text{Var}\!\left(#1\right)}
\newcommand{\VAvarover}[2]{\text{Var}_{#1}\!\left(#2\right)}
\newcommand{\VAcov}[1]{\text{Cov}\!\left(#1\right)}

\definecolor{algcomment}{gray}{0.25}

\algrenewcommand\alglinenumber[1]{\footnotesize #1}
\algrenewcommand{\algorithmiccomment}[1]{{\color{algcomment}\hfill$\rightarrow$ #1}}

\newcommand{\VAaprealoperator}{\mathcal{S}}
\newcommand{\VAaprealoperatorD}{\VAaprealoperator^{>0}}

\newcommand{\VAapkindifoperator}{\tilde{\mathcal{S}}}
\newcommand{\VAapDoperator}{\widehat{\mathcal{S}}}
\newcommand{\VAapKoperatoronlyD}{{\mathcal{S}}'}

\newcommand{\VAapKDoperator}{\VAapkindifoperator}
\newcommand{\VAapKDoperatorD}{\VAapKDoperator^{>0}}

\newcommand{\VAapkinkinoperatorD}{\VAaprealoperator^{=1}}

\newcommand{\zelfkalk}[1]{\mathcal{W}_1\left(#1\right)}
\newcommand{\zelfkalkdun}[1]{\mathcal{W}_1\!\left(\!#1\!\right)\!}

\newcommand{\VAfinalV}{\nu}

\makeatletter
\algrenewcommand\ALG@beginalgorithmic{\normalsize}
\algrenewcommand\algorithmiccomment[2][\footnotesize]{{#1\hfill\(\triangleright\) #2}}
\makeatother

\begin{document}
\maketitle

\section*{Abstract}

We develop a novel Monte Carlo strategy for the simulation of the Boltzmann-BGK model with both low-collisional and high-collisional regimes present. The presented solution to maintain accuracy in low-collisional regimes and remove exploding simulation costs in high-collisional regimes uses hybridized particles that exhibit both kinetic behaviour and diffusive behaviour depending on the local collisionality. In this work, we develop such a method that maintains the correct mean, variance, and correlation of the positional increments over multiple time steps of fixed step size for all values of the collisionality, under the condition of spatial homogeneity during the time step. In the low-collisional regime, the method reverts to the standard velocity-jump process. In the high-collisional regime, the method collapses to a standard random walk process. We analyze the error of the presented scheme in the low-collisional regime for which we obtain the order of convergence in the time step size. We furthermore provide an analysis in the high-collisional regime that demonstrates the asymptotic-preserving property.

Keywords: Monte Carlo, asymptotic preserving, Boltzmann-BGK, kinetic-diffusion

\section{Introduction}

Applications such as rarefied gas modeling~\cite{boyddeschenes2011aerospace}, radiation transport~\cite{fleck1971implMCphotontrans}, and neutral or ion transport in nuclear fusion~\cite{larsen1974neutrontransportbecomesdiff,stangeby2000plasmaboundary} frequently have to cope with large differences in collision rates between simulated regions. In low-collisional regions, a kinetic description is often required for accuracy, whereas this kinetic description becomes computationally intractable in high-collisional regions. However, in high-collisional regimes, a limiting fluid description can become valid, which is cheaper to simulate. One strategy to handle such situations uses domain decomposition, where part of the domain is described via the kinetic model and part by the fluid description~\cite{boyddeschenes2011aerospace,densmore2007domdecMC}. Another popular method is separating the density in both a kinetic part and a fluid part throughout the domain~\cite{pareschi1999implicitMCrarefied,crouseilles2004hybridgas,dimarco2008hybridIIkin,horsten2018hybrid}. Both types of solutions involve complications in determining the partition of the domain and in coupling the two different parts of the model.


A different solution type is formed by the asymptotic-preserving (AP) methods that avoid couplings by using a single method throughout the domain. Such methods are designed with the accuracy of the kinetic simulation in the kinetic regions and the efficiency of a fluid simulation in the fluid regions. The first such method was developed in the context of radiation transport~\cite{fleck1971implMCphotontrans,fleck1984APpBrownian}, for neutron transport~\cite{borgers1992asymptoticdiffusionLTE}, and for the Boltzmann equation~\cite{gabetta1997timerelaxation}. Most of these methods are fully deterministic~\cite{bennoune2008apboltmzannNS,boscarino2013apdiflimRK,buet2007apradtransfer,
crouseilles2011apmMvlasov,dimarco2012highorderapboltzmann,klar1998asinddifflim,
jin1998diffrelaxdiscrvel,lemou2008APTLEdiffmM,naldi2000aphyperbolicdifflim} 
and generally fully resolve the velocity domain, which is unnecesary in the fluid limit. Asymptotic preserving Monte Carlo methods (APMC)~\cite{pareschi2001timerelaxed,gorji2014particleFPrarefied,dimarco2018APdiff,bufferand2013plasmaparallelheat,dicintio2017plasmaparallelheatextrafactor,crestetto2018APdiff} avoid unnecessarily resolving the velocity dimensions and have as additional advantages dimension independence and the capability to easily cope with complex geometries.

Here, we develop a new kinetic-diffusion APMC method that combines a standard Monte Carlo (MC) method for the kinetic equation with a random walk MC method for the limiting advection-diffusion equation. In each time step, the MC particles in the new method move according to the kinetic equation until they collide. After a collision, the MC particle moves diffusively according to a random walk process for the remainder of that time step. This diffusive motion is such that the positional increments match the kinetic process' mean and variance exactly in every time step. Furthermore, the correlation of the motion between subsequent time steps is maintained by the proposed combination of kinetic and advection-diffusion parts.

In the low-collisional regime, the kinetic parts of the methods prevail, resulting in only very minor parts of the particle paths where the particle moves according to the (in that limit invalid) diffusive process. In the high-collisional regime, the kinetic parts make up only a marginal part of the process, resulting in a large gain in efficiency. The presented method is furthermore easily implementable due to its simple structure: kinetic motion with filling of the time steps with a diffusive step.

In Section~\ref{sec:ap_kinetic}, we present the kinetic description of the Boltzmann-BGK model and the corresponding standard MC simulation, which forms the basis of the new algorithm. In Section~\ref{sec:ap_aggregtodiff}, we derive expressions for the mean and variance of a positional increment in the standard MC simulation of the kinetic model. These expressions are used for the advection-diffusion coefficients in the diffusive motion in the kinetic-diffusion (KD) simulation method. The KD simulation method is presented in Section~\ref{sec:ap_newscheme}. In Sections~\ref{sec:ap_limnul} and~\ref{sec:ap_liminf}, we provide error bounds on the presented simulation schemes in respectively the low-collisional and high-collisional limit, proving both consistency and the asymptotic-preserving property. Finally, in Section~\ref{sec:ap_num}, a numerical illustration of the low error and the speed-up is presented.

\section{The kinetic model and its simulation\label{sec:ap_kinetic}}

In Section~\ref{subsec:ap_kinetic_kinetic}, we present the kinetic model in its integro-differential form. In Section~\ref{subsec:ap_kinetic_standardMC}, we introduce the standard Monte Carlo method for the kinetic model, which is founded in the particle description of the kinetic model. Then, in Section~\ref{subsec:ap_modelsim_difflimit}, we consider the limiting equation in the diffusive limit. Finally, in Section~\ref{subsec:ap_modelsim_strategy}, we discuss the exploding simulation cost of the standard Monte Carlo method in diffusive regimes and shortly introduce our strategy to mitigate this computational burden, which will be based on the limiting behaviour as found in Section~\ref{subsec:ap_modelsim_difflimit}.

\def\VApostcolveldistri{\mathcal{M}}
\def\VAbackgroundspeedtemp{T}
\def\VAbackgroundspeedmeanflow{u}
\def\VAratespatial{\sigma}

\subsection{Kinetic model\label{subsec:ap_kinetic_kinetic}}

The BGK kinetic model describes particles that move according to a velocity jump process: the particle moves with a constant velocity until it collides, at which point its velocity changes to a new velocity, sampled from the position-dependent distribution $\VApostcolveldistri(v;x)$. We write the integro-differential equation for the density $f(x,v,t)$ of such particles as
\begin{align}
\frac{\partial f(x,v,t)}{\partial t}+\underbrace{\frac{v}{\VAscaleparam}\frac{\partial f(x,v,t)}{\partial x}\vphantom{\left(\int\right)}}_\text{transport}&=\underbrace{\frac{\VAratespatial(x)}{\VAscaleparam^2}\left(\VApostcolveldistri(v;x)\int f(x,v,t)\text{d}v-f(x,v,t)\right)}_\text{collision}\,,\label{eq:ap_kinetic_integrodiff}\\
f(x,v,0)&=\VAnsource(x,v)\,,
\end{align}
with initial condition $\VAnsource(x,v)$, scattering rate $\frac{\VAratespatial(x)}{\VAscaleparam^2}$, scaling parameter $\VAscaleparam$, and post-collisional velocity distribution $\VApostcolveldistri(v;x)$. The scaling parameter $\VAscaleparam$ captures the diffusive scaling~\cite{bardos1993kintofluid} as $\VAscaleparam\rightarrow0$, which describes a situation with large velocities and a very high collision rate. The post-collisional velocity distribution, $\VApostcolveldistri(v;x)$, is often a Maxwellian and we will denote the mean velocity as $\VAbackgroundspeedmeanflow(x)$ and the temperature as $\VAbackgroundspeedtemp(x)$, leading to a probability density function
\begin{equation}
\VApostcolveldistri(v;x)=\frac{1}{\sqrt{2\pi\VAbackgroundspeedtemp(x)}}e^{-\frac{(v-\VAscaleparam\VAbackgroundspeedmeanflow(x))^2}{2\VAbackgroundspeedtemp(x)}}\,.\label{eq:ap_kinetic_postcolveldistr_pdf}
\end{equation}
The post-collisional velocity distribution being normal is not a requirement for the method proposed in this manuscript, but knowledge of the mean plasma velocity $\VAbackgroundspeedmeanflow(x)$ and temperature $\VAbackgroundspeedtemp(x)$ is essential.

In a general 3D setting, the dimensionality of this equation is seven, prompting the use of Monte Carlo methods. We restrict our discussion to the 1D version of the equation, since both the standard Monte Carlo method and the simulation method we propose can be readily extended to higher dimensions by treating the different dimensions identically and independently.

\subsection{Kinetic particle description and standard Monte Carlo\label{subsec:ap_kinetic_standardMC}}

Underlying Equation~\eqref{eq:ap_kinetic_integrodiff} is a particle description which can be used for a standard MC method and which is the topic of this section.

We denote the discrete times at which the particle velocity changes as $\VAeventtime_\VAeventno$, $\VAeventno\in\{1,2,\dots\}$. The equations for the position $\VAtimevar\mapsto {x}(\VAtimevar)$ and for the velocity $\VAtimevar\mapsto \frac{{v}(\VAtimevar)}{\VAscaleparam}$ can be written as 
\begin{align}
\left({x}(0),{v}(0)\right)&\sim \VAnsource({x},{v})\\
\frac{\text{d}{x}(\VAtimevar)}{\text{d}t}&=\frac{{v}(\VAtimevar)}{\VAscaleparam}\\
{v}(\VAtimevar)&=\VAdiscreteveleend_\VAeventno\text{ for }\VAtimevar\in\left[\VAeventtime_{\VAeventno},\VAeventtime_{\VAeventno+1}\right),\ k\in\left\{0,\dots,\VAnofevents-1\right\}\,.
\end{align}



The collisions with the plasma occur with position-dependent rate $\frac{\VAratespatial({x})}{\VAscaleparam^2}$. 
Sampling the event times can be done by using standard exponentially distributed samples $\VAstandardexp_\VAeventno\sim\mathcal{E}(1)$ and solving the equation
\begin{equation}
\int_{0}^{\VAeventtimedif_{\VAeventno}}\frac{\VAratespatial\left({x}(\VAeventtime_\VAeventno)+\frac{\VAdiscreteveleend_\VAeventno}{\VAscaleparam}\VAtimevar\right)}{\VAscaleparam^2}\text{d}\VAtimevar=\VAstandardexp_\VAeventno\,,\label{eq:ap_kin_eventtime}
\end{equation}
for the free-flight intervals $\VAeventtimedif_{\VAeventno}$, from which the next event time can be found as $\VAeventtime_{\VAeventno+1}=\VAeventtime_\VAeventno+\VAeventtimedif_\VAeventno$. We will refer to this sampling of $\VAstandardexp_\VAeventno$ and the following inversion of the integral in Equation \eqref{eq:ap_kin_eventtime} by the function \textproc{SampleCollision}$({x}(\VAeventtime_\VAeventno),\VAdiscreteveleend_\VAeventno,\VAratespatial(x),\VAscaleparam)$.\label{plaats:ap_kin_eventtime_function} 
At a collision, a new velocity is sampled according to the post-collisional velocity distribution, which is a Maxwellian distribution, Equation~\eqref{eq:ap_kinetic_postcolveldistr_pdf}, in this text:
\begin{equation}
\VAdiscreteveleend_\VAeventno\sim\VApostcolveldistri(v;x(\VAeventtime_\VAeventno))\,.
\end{equation}

A Monte Carlo algorithm to simulate Equation~\eqref{eq:ap_kinetic_integrodiff} based on this particle description of the kinetic model up to the final time $\VAtotaltime$ is algorithm~\ref{alg:ap_kin}.

\algrenewcommand\alglinenumber[1]{\footnotesize #1}
\begin{algorithm}
  \caption{A kinetic simulation up to time $\VAtotaltime$
    \label{alg:ap_kin}}
  \begin{algorithmic}[1]
    \Function{KineticSimulation}{$\VAdiscreteposeend,\VAdiscreteveleend,\VAtotaltime,\VAratespatial({x}),\VAscaleparam,\VApostcolveldistri(v;x)$}
      \State $\VAtimevar\gets0$
      \While{$\VAtimevar<\VAtotaltime$}
		\State $\VAeventtimedif\gets$ \Call{SampleCollision}{$\VAdiscreteposeend,\VAdiscreteveleend,\VAratespatial({x}),\VAscaleparam$}
		\State $\tau\gets$\Call{min}{$\VAeventtimedif,\VAtotaltime-\VAtimevar$}\Comment{determine the kinetically moved time}
		\State $\VAdiscreteposeend\gets \VAdiscreteposeend+\frac{\VAdiscreteveleend}{\VAscaleparam}\tau$\Comment{execute the free flight}
        \If{$\VAeventtimedif<\VAtotaltime-\VAtimevar$}\Comment{check if a collision occurred}
          \State $\VAdiscreteveleend\gets\VApostcolveldistri(v;x)$\Comment{sample the post-collisional velocity}
        \EndIf
		\State $\VAtimevar\gets\VAtimevar+\tau$
      \EndWhile
      \State \Return{$\VAdiscreteposeend,\VAdiscreteveleend$}
    \EndFunction
  \end{algorithmic}
\end{algorithm}

\newcommand{\VAdensity}{\rho}

\subsection{Advection-diffusion limit in the high-collisional regime\label{subsec:ap_modelsim_difflimit}}

In this section, we will uncover the limiting behaviour of Equation~\eqref{eq:ap_kinetic_integrodiff} when $\VAscaleparam\rightarrow0$, according to a standard strategy such as presented in~\cite{othmer2000veljumpdifflimit}.

When the scaling parameter $\VAscaleparam$ becomes small, the collision rate $\frac{\VAratespatial(x)}{\VAscaleparam^2}$ will become very large, resulting in a very quick relaxation of the distribution $f(x,v,t)$ to $\VApostcolveldistri(v;x)\VAdensity(x,t)$ where $\VAdensity(x,t)=\int f(x,v,t)\text{d}v$ represents the density. We capture this behaviour by writing the distribution $f(x,v,t)$ as
\begin{equation}
f(x,v,t)=\VApostcolveldistri(v;x)\VAdensity(x,t)+\VAscaleparam g(x,v,t)\,,\label{eq:ap_modelsim_difflimit_pluscorr}
\end{equation}
i.e., the distribution $\VApostcolveldistri(v;x)\VAdensity(x,t)$ to which $f(x,v,t)$ relaxes plus a small correction term $\VAscaleparam g(x,v,t)$. Plugging Equation~\eqref{eq:ap_modelsim_difflimit_pluscorr} in Equation~\eqref{eq:ap_kinetic_integrodiff} gives
\begin{equation}
\frac{\partial\!\left(\!\VApostcolveldistri(v;x)\!\VAdensity(x,t)\!+\!\VAscaleparam g(x,v,t)\!\right)\!}{\partial t}+\frac{v}{\VAscaleparam}\frac{\partial\!\left(\!\VApostcolveldistri(v;x)\!\VAdensity(x,t)\!+\!\VAscaleparam g(x,v,t)\!\right)\!}{\partial x}=-\frac{\VAratespatial(x)}{\VAscaleparam}g(x,v,t)\,,\label{eq:ap_modelsim_difflimit_pluggedin}
\end{equation}
which, after averaging over $v$ becomes
\begin{equation}
\frac{\partial\VAdensity(x,t)}{\partial t}+\frac{1}{\VAscaleparam}\frac{\partial\left(\int v\VApostcolveldistri(v;x)\text{d}v\VAdensity(x,t)\right)}{\partial x}+\int v\frac{\partial g(x,v,t)}{\partial x}\text{d}v=0\,.\label{eq:ap_modelsim_difflimit_averaged}
\end{equation}
When we consider the dominant terms of Equation~\eqref{eq:ap_modelsim_difflimit_pluggedin}, which have $\frac{1}{\VAscaleparam}$ as a factor, we find that
\begin{equation}
v\frac{\partial (\VApostcolveldistri(v;x)\VAdensity(x,t))}{\partial x}=-\VAratespatial(x) g(x,v,t)\,,
\end{equation}
which can be used to transform Equation~\eqref{eq:ap_modelsim_difflimit_averaged} into
\begin{equation}
\frac{\partial\VAdensity(x,t)}{\partial t}+\frac{\partial\left(\VAbackgroundspeedmeanflow(x)\VAdensity(x,t)\right)}{\partial x}-\frac{\partial}{\partial x}\left(\frac{1}{\VAratespatial(x)}\frac{\partial(\VAbackgroundspeedtemp(x)\VAdensity(x,t))}{\partial x}\right)=0\,,\label{eq:ap_modelsim_difflimit_advdif}
\end{equation}
as $\VAscaleparam\rightarrow0$. Equation~\eqref{eq:ap_modelsim_difflimit_advdif} is an advection-diffusion-type equation, which validates the name {diffusive limit}, when $\VAscaleparam\rightarrow0$. Note that the effects of the collision rate, $\VAratespatial(x)$, and the temperature, $\VAbackgroundspeedtemp(x)$, differ in nature, since these coefficients appear at different places in the equation.

In a spatially homogeneous setting with $\VAbackgroundspeedmeanflow(x)\equiv\VAbackgroundspeedmeanflow$, $\VAbackgroundspeedtemp(x)\equiv\VAbackgroundspeedtemp$, and $\VAratespatial(x)\equiv\VAratespatial$, Equation~\eqref{eq:ap_modelsim_difflimit_advdif} is the Fokker-Planck equation of the stochastic differential equation (SDE)
\begin{equation}
\text{d}X=u\text{d}t+\sqrt{\frac{2\VAbackgroundspeedtemp}{\VAratespatial}}\text{d}W\,,\label{eq:ap_modelsim_difflimit_sdeequiv}
\end{equation}
with $\text{d}W$ expressing a standard Wiener process.

\subsection{Exploding simulation costs in the diffusive limit\label{subsec:ap_modelsim_strategy}}

When we consider the kinetic simulation presented in Section~\ref{subsec:ap_kinetic_standardMC} in the limit $\VAscaleparam\rightarrow0$, the collision rate $\frac{\VAratespatial}{\VAscaleparam^2}$ becomes very large and therefore the execution of all individual collisions and free flight becomes computationally expensive. The simulation of the SDE of Equation~\eqref{eq:ap_modelsim_difflimit_sdeequiv}, however, can be done at a computational cost that does not depend on $\VAscaleparam$ via a random walk simulation. This random walk simulation of Equation~\eqref{eq:ap_modelsim_difflimit_sdeequiv} with a fixed time step $\VAtimestep$ amounts to
\begin{equation}
\VAdiscreteposeend(t+\VAtimestep)=\VAdiscreteposeend(t)+A_0\VAtimestep+\sqrt{2D_0\VAtimestep}\xi\,,\label{eq:ap_modelsim_difflimit_sderandomwalk}
\end{equation}
with $A_0=\VAbackgroundspeedmeanflow$, $D_0=\frac{\VAbackgroundspeedtemp}{\VAratespatial}$, and $\xi\sim\mathcal{N}(0,1)$. However, use of Equations~\eqref{eq:ap_modelsim_difflimit_sdeequiv} and~\eqref{eq:ap_modelsim_difflimit_sderandomwalk} is only allowed in the limit $\VAscaleparam\rightarrow0$, and therefore introduces a modeling error for finite $\VAscaleparam$.

To cope with the high computational cost for small but finite $\VAscaleparam$, we aim at replacing most free flights in the high-collisional regime by a single random walk step with $\VAscaleparam$-dependent coefficients $A_\VAscaleparam$ and $D_\VAscaleparam$ that covers a fraction of the time step. This corresponds to applying a random walk Monte Carlo discretization of an advection-diffusion equation for that fraction of the time step. To remove the modeling error for finite $\VAscaleparam$ as the time step $\VAtimestep$ decreases, we propose two kinetic corrections in the algorithm of Section~\ref{sec:ap_newscheme}.

\newcommand{\VAdiftime}{\theta}
The first correction to control the modeling error is not executing the entirety of the time step $\VAtimestep$ diffusively, but maintaining kinetic motion at the beginning of each time step. The resulting algorithm presented in Section~\ref{sec:ap_newscheme} is consequently called the kinetic-diffusion (KD) scheme. As a second correction, we do not use the limiting advection and diffusion coefficients $\VAbackgroundspeedmeanflow$ and $\VAbackgroundspeedtemp$, but instead use coefficients $A_\VAscaleparam$ and $D_\VAscaleparam$ such that the first two moments of the diffusive motion match the kinetic process exactly for any fixed and finite value of $\VAscaleparam$. Then, the diffusive substep of the new Monte Carlo algorithm becomes
\begin{equation}
\VAdiscreteposeend(t+\VAdiftime)=\VAdiscreteposeend(t)+A_\VAscaleparam\VAdiftime+\sqrt{2D_\VAscaleparam\VAdiftime}\xi\,,\label{eq:ap_modelsim_strategy_randomwalk}
\end{equation}
with $0\leq\VAdiftime\leq\VAtimestep$ the diffusive part of the stepped time and $\xi\sim\mathcal{N}(0,1)$ a standard normal distributed random number. The expressions for the advection and diffusion coefficients $A_\VAscaleparam$ and $D_\VAscaleparam$ are derived in Section~\ref{sec:ap_aggregtodiff}. In Section~\ref{sec:ap_newscheme}, the new kinetic-diffusion scheme is presented, which is analyzed in the remainder of the paper.

\section{Mean and variance of the kinetic motion\label{sec:ap_aggregtodiff}}

In this Section, we will derive the mean and variance of the motion of a kinetically simulated particle during a time step of length $\VAtimestep$ for a spatially homogeneous background. In Section~\ref{sec:ap_newscheme}, these moments will be used for the advection and diffusion coefficients $A_\VAscaleparam$ and $D_\VAscaleparam$ for positional increments. The assumption of a spatially homogeneous background translates to a constant rate $\frac{\VAratespatial}{\VAscaleparam^2}$, and a post-collisional velocity distribution with constant mean $\VAscaleparam\VAbackgroundspeedmeanflow$ and variance $\VAbackgroundspeedtemp$. 

\begin{figure}[H]
\centering
\begin{minipage}{\textwidth}
  \centering
  \resizebox{\textwidth}{!}{%
	\tikzset{onderaannode/.style={font=\small}}
\begin{tikzpicture}[outer sep=0pt]
    \begin{axis}[
     axis lines=left,
        xmin=-2.5,xmax=25.5,
  ticks=none,
  xlabel={},
        ylabel={},
        ymin=-1,ymax=1,
        width=1.3*\textwidth,
        height=.3\textwidth,
        axis line style={draw=none},
        legend style={legend cell align=right,legend plot pos=right}]

    \def\zelfH{.2}
    \def\zelfAH{.2}
    \def\zelfTHD{.45}
    \def\zelfAS{.05}
 \draw[-{Latex[length=.3cm,width=.2cm]}] (axis cs:-3,0) to (axis cs:25.5,0);
 \node at (axis cs:25.2,.3) {$t$};
 \draw (axis cs:2,\zelfH) to (axis cs:2,-\zelfH);
 \node at (axis cs:2,\zelfTHD) {$\VAtimestepno\VAtimestep$};
 \draw (axis cs:18,\zelfH) to (axis cs:18,-\zelfH);
 \node at (axis cs:18,\zelfTHD) {$(\VAtimestepno+1)\VAtimestep$};
 \foreach \i/\j in {1/5,2/8}{
  \edef\temp{\noexpand\draw plot[mark=*, mark options={color=black!20!red}, mark size=3] coordinates {(axis cs:{\j},0)};}
  \temp
  \edef\temp{\noexpand\node at (axis cs:{\j},\zelfTHD) {$\VAeventtime_{\VAeventnotilltimestep{\VAtimestepno}+\i}$};}
     \temp}
    \draw plot[mark=*, mark options={color=black!20!red}, mark size=3] coordinates {(axis cs:{12},0)};
    \node at (axis cs:{-2},\zelfTHD) {$\VAeventtime_{\VAeventnotilltimestep{\VAtimestepno}}$};
    \node at (axis cs:{12},\zelfTHD) {$\VAeventtime_{\VAeventnotilltimestep{\VAtimestepno}+3}=\VAeventtime_{\VAeventnotilltimestep{{\VAtimestepno+1}}}$};
    \node at (axis cs:{23.5},\zelfTHD) {$\VAeventtime_{\VAeventnotilltimestep{\VAtimestepno+1}+1}$};
 \draw plot[mark=*, mark options={color=black}, mark size=3] coordinates {(axis cs:-2,0)};
 \draw plot[mark=*, mark options={color=black}, mark size=3] coordinates {(axis cs:23.5,0)};
 \path let \p1=($(axis cs:0,.2)-(axis cs:0,0)$) in \pgfextra{\xdef\myyshift{\y1}};
 \draw [decorate,decoration={brace,amplitude=5pt},xshift=0,yshift=-\myyshift/2] (axis cs:{2-\zelfAS},-\zelfAH) -- (axis cs:{-2+\zelfAS},-\zelfAH) node [onderaannode,black,midway,yshift=-\myyshift*2.5] {$\VAeventtimedif_{\VAeventnotilltimestep{\VAtimestepno}+1|n-1}$};
 \draw [decorate,decoration={brace,amplitude=5pt},xshift=0,yshift=-\myyshift/2] (axis cs:{5-\zelfAS},-\zelfAH) -- (axis cs:{2+\zelfAS},-\zelfAH) node [onderaannode,black,midway,yshift=-\myyshift*2.5] {$\VAeventtimedif_{\VAeventnotilltimestep{\VAtimestepno}+1|n}$};
 \draw [decorate,decoration={brace,amplitude=5pt},xshift=0,yshift=-\myyshift/2] (axis cs:{8-\zelfAS},-\zelfAH) -- (axis cs:{5+\zelfAS},-\zelfAH) node [onderaannode,black,midway,yshift=-\myyshift*2.5] {$\VAeventtimedif_{\VAeventnotilltimestep{\VAtimestepno}+2}$};
 \draw [decorate,decoration={brace,amplitude=5pt},xshift=0,yshift=-\myyshift/2] (axis cs:{12-\zelfAS},-\zelfAH) -- (axis cs:{8+\zelfAS},-\zelfAH) node [onderaannode,black,midway,yshift=-\myyshift*2.5] {$\VAeventtimedif_{\VAeventnotilltimestep{\VAtimestepno}+3}$};
 \draw [decorate,decoration={brace,amplitude=5pt},xshift=0,yshift=-\myyshift/2] (axis cs:{18-\zelfAS},-\zelfAH) -- (axis cs:{12+\zelfAS},-\zelfAH) node [onderaannode,black,midway,yshift=-\myyshift*2.5] {$\ \ \ \VAeventtimedif_{\VAeventnotilltimestep{\VAtimestepno}+4|n}=\VAeventtimedif_{\VAeventnotilltimestep{\VAtimestepno+1}+1|n}$};
 \draw [decorate,decoration={brace,amplitude=5pt},xshift=0,yshift=-\myyshift/2] (axis cs:{23.5-\zelfAS},-\zelfAH) -- (axis cs:{18+\zelfAS},-\zelfAH) node [onderaannode,black,midway,yshift=-\myyshift*2.5] {$\ \ \ \VAeventtimedif_{\VAeventnotilltimestep{\VAtimestepno+1}+1|n+1}$};
    \end{axis}%
\end{tikzpicture}%
}
  \captionof{figure}{The part of the free flight times (time between collisions) that overlap with the $\VAtimestepno$-th time step for an example with three collisions in the time interval $[\VAtimestepno\VAtimestep,(\VAtimestepno+1)\VAtimestep)$.}
  \label{fig:ap_eventtime_tijdsstap}
\end{minipage}\qquad
\end{figure}
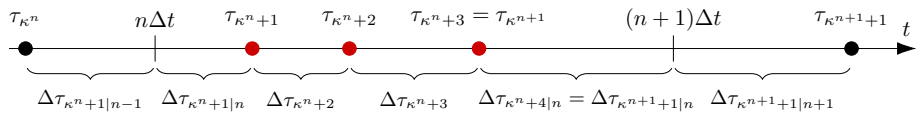

During the $\VAtimestepno$-th time step $[\VAtimestepno\VAtimestep,(\VAtimestepno+1)\VAtimestep)$ the change in position of a particle is the result of the free flights during that time step. In Figure~\ref{fig:ap_eventtime_tijdsstap}, these free flights occur in the time between the collisions, which are indicated with circles. Only the parts of the free flights that overlap with the time step contribute to the motion in the time step. We denote the index of the last collision before the time interval $[\VAtimestepno\VAtimestep,(\VAtimestepno+1)\VAtimestep)$ with $\VAeventnotilltimestep{\VAtimestepno}$ and we denote the part of the $\VAeventnotilltimestep{\VAtimestepno}+1$-th flight time that is in the time interval $[\VAtimestepno\VAtimestep,(\VAtimestepno+1)\VAtimestep)$ with $\VAeventtimedif_{\VAeventnotilltimestep{\VAtimestepno}+1|\VAtimestepno}$ and similarly for the overlap of the $\VAeventnotilltimestep{\VAtimestepno}+1$-th flight time with the time interval $[\VAtimestepno\VAtimestep,(\VAtimestepno+1)\VAtimestep)$, we write $\VAeventtimedif_{\VAeventnotilltimestep{\VAtimestepno}+1|\VAtimestepno}$. With this notation, the change in position during the $\VAtimestepno$-th time step can be written as
{
\newcommand{\zelflocalmaxvertical}{\vphantom{\sum_{\VAeventnointimestep=\VAeventnotilltimestep{\VAtimestepno}+2}^{\VAeventnotilltimestep{\VAtimestepno+1}}}}
\begin{equation}
\Delta {\VAdiscreteposeend}_\VAtimestepno=\underbrace{\zelflocalmaxvertical\frac{\VAdiscreteveleend_{\VAeventnotilltimestep{\VAtimestepno}}}{\VAscaleparam}\VAeventtimedif_{\VAeventnotilltimestep{\VAtimestepno}+1|\VAtimestepno}}_\text{entering flight}+\!\!\underbrace{\zelflocalmaxvertical\sum_{\VAeventnointimestep=\VAeventnotilltimestep{\VAtimestepno}+2}^{\VAeventnotilltimestep{\VAtimestepno+1}}\!\!\frac{\VAdiscreteveleend_{\VAeventnointimestep-1}}{\VAscaleparam}\VAeventtimedif_\VAeventnointimestep}_\text{internal flights}+\underbrace{\zelflocalmaxvertical\frac{\VAdiscreteveleend_{\VAeventnotilltimestep{\VAtimestepno+1}}}{\VAscaleparam}\VAeventtimedif_{\VAeventnotilltimestep{\VAtimestepno+1}+1|n}}_\text{exiting flight}\label{eq:ap_deltax_real}
\end{equation}
}It is important to note that in this model, there is correlation between the positional increments during different time steps, because the last velocity of a time step, $\VAdiscreteveleend_{\VAeventnotilltimestep{\VAtimestepno}}$, is also the first velocity of the next time step.

Since we are only considering a single time step, it is convenient to discard the dependence on the time step index, $\VAtimestepno$. We rename the velocities $\{\VAdiscreteveleend_{\VAeventnotilltimestep{\VAtimestepno}},\cdots,\VAdiscreteveleend_{\VAeventnotilltimestep{\VAtimestepno+1}}\}$ as $\{\VAdiscreteveleend_0,\cdots,\VAdiscreteveleend_\VAnofeventsintimestepnon\}$ and the flight times $\{\VAeventtimedif_{\VAeventnotilltimestep{\VAtimestepno}+1|n},\VAeventtimedif_{\VAeventnotilltimestep{\VAtimestepno}+2},\cdots,\VAeventtimedif_{\VAeventnotilltimestep{\VAtimestepno+1}},\VAeventtimedif_{\VAeventnotilltimestep{\VAtimestepno+1}+1|n}\}$ as $\{\VAeventtimedif_0,\cdots,\VAeventtimedif_{\VAnofeventsintimestepnon}\}$, with $\VAnofeventsintimestepnon=\VAeventnotilltimestep{\VAtimestepno+1}-\VAeventnotilltimestep{\VAtimestepno}$. This rephrases equation~\eqref{eq:ap_deltax_real} as
\begin{equation}
\Delta \VAdiscreteposeend=\sum_{\VAeventnointimestep=0}^\VAnofeventsintimestepnon\frac{\VAdiscreteveleend_\VAeventnointimestep}{\VAscaleparam}\VAeventtimedif_\VAeventnointimestep\,,\label{eq:ap_deltax_assum}
\end{equation}
where we always have
\begin{equation}
\sum_{\VAeventnointimestep=0}^\VAnofeventsintimestepnon\VAeventtimedif_\VAeventnointimestep=\VAtimestep\,.\label{eq:ap_deltat_assum}
\end{equation}
Note that $\VAnofeventsintimestepnon$, the number of collisions during a time step of length $\VAtimestep$, is a random variable which follows the Poisson distribution:
\begin{equation}
P(\VAnofeventsintimestepnon=\VAeventno)=\frac{1}{\VAeventno!}\left(\frac{\VAratespatial}{\VAscaleparam^2}\VAtimestep\right)^{\VAeventno}e^{-\frac{\VAratespatial}{\VAscaleparam^2}\VAtimestep}\,.\label{eq:ap_poissondistri}
\end{equation}

We will now derive the mean and variance of $\Delta \VAdiscreteposeend$. In Section~\ref{subsec:ap_randomv0}, we first neglect correlation with previous and next time steps by treating the initial and final velocity identically to the other velocities. From the resulting formulas from Section~\ref{subsec:ap_randomv0}, we proceed to derive formulas for the mean and variance of the positional increment conditioned on a known final velocity, $\VAdiscreteveleend_\VAnofeventsintimestepnon$, in Section~\ref{subsec:ap_fixedv0}. The formulas resulting from this second section allow incorporating correlation with subsequent time steps in the random walk steps itself, leading to the algorithm of Section~\ref{sec:ap_newscheme}.

\subsection{Neglecting correlation with other time steps\label{subsec:ap_randomv0}}

In this first section, we derive the exact mean and variance of the positional increment of a particle adhering to the model of Section~\ref{subsec:ap_kinetic_kinetic} during a time step, while neglecting correlation between different time steps. The resulting expressions form the foundation of the derivation in Section~\ref{subsec:ap_fixedv0}, where we do consider correlation with the subsequent time step.

\subsubsection{Neglecting correlation with other time steps: mean\label{subsubsec:ap_randomv0_mean}}

The change in position of a particle following the model of Section~\ref{subsec:ap_kinetic_kinetic} during a time step $\VAtimestep$ is written as $\Delta\VAdiscreteposeend$ and can be expressed as in Equation~\eqref{eq:ap_deltax_assum}. With this description, and by conditioning on the Poisson-distributed number of collisions $\VAnofeventsintimestepnon$, we can write the expected value of $\Delta\VAdiscreteposeend$ as
\begin{equation}
\VAexpec{\Delta \VAdiscreteposeend}=\VAexpec{\sum_{\VAeventnointimestep=0}^\VAnofeventsintimestepnon\frac{\VAdiscreteveleend_\VAeventnointimestep}{\VAscaleparam}\VAeventtimedif_\VAeventnointimestep}=\VAexpecover{\VAnofeventsintimestepnon}{\VAexpec{\left.\sum_{\VAeventnointimestep=0}^\VAnofeventsintimestepnon\frac{\VAdiscreteveleend_\VAeventnointimestep}{\VAscaleparam}\VAeventtimedif_\VAeventnointimestep\right|\VAnofeventsintimestepnon}}\,,\label{eq:ap_randomv0_mean_initial}
\end{equation}
where $\mathbb{E}_{\VAnofeventsintimestepnon}$ is the expected value over the specific random variable in the subscript and the expected value without subscript goes over all random variables except the ones on which the expected value is conditioned.

The different $\VAdiscreteveleend_\VAeventno$ in Equation~\eqref{eq:ap_randomv0_mean_initial} are random variables that follow the distribution $\VApostcolveldistri(v)$ with mean $\VAscaleparam\VAbackgroundspeedmeanflow$ and variance $\VAbackgroundspeedtemp$ and the different $\VAdiscreteveleend_\VAeventno$ are independent of all other random variables. Using this transforms Equation~\eqref{eq:ap_randomv0_mean_initial} into
\begin{equation}
\VAexpec{\left.\sum_{\VAeventnointimestep=0}^\VAnofeventsintimestepnon\frac{\VAdiscreteveleend_\VAeventnointimestep}{\VAscaleparam}\VAeventtimedif_\VAeventnointimestep\right|\VAnofeventsintimestepnon}=\VAexpec{\left.\sum_{\VAeventnointimestep=0}^\VAnofeventsintimestepnon\VAexpec{\frac{\VAdiscreteveleend_\VAeventnointimestep}{\VAscaleparam}}\VAeventtimedif_\VAeventnointimestep\right|\VAnofeventsintimestepnon}=\VAbackgroundspeedmeanflow\VAexpec{\left.\sum_{\VAeventnointimestep=0}^\VAnofeventsintimestepnon\VAeventtimedif_\VAeventnointimestep\right|\VAnofeventsintimestepnon}\,.\label{eq:ap_randomv0_mean_vdistriworkedout}
\end{equation}
Since the $\VAeventtimedif_\VAeventnointimestep$ always sum to $\VAtimestep$ by construction of $\Delta \VAdiscreteposeend$ (see Equation~\eqref{eq:ap_deltat_assum}) regardless of $\VAnofeventsintimestepnon$, we find
\begin{equation}
\VAexpec{\Delta \VAdiscreteposeend}=\VAbackgroundspeedmeanflow\VAtimestep\,.\label{eq:ap_randomv0_mean_result}
\end{equation}

\subsubsection{Neglecting correlation with other time steps: variance\label{subsubsec:ap_fixedv0_var}}

To find an expression for the variance of the positional increment, we start from the law of total variance,
\begin{equation}
\VAvar{\Delta \VAdiscreteposeend}=\VAexpecover{\VAnofeventsintimestepnon}{\VAvar{\left.\sum_{\VAeventnointimestep=0}^\VAnofeventsintimestepnon \frac{\VAdiscreteveleend_\VAeventnointimestep}{\VAscaleparam}\VAeventtimedif_\VAeventnointimestep\right|\VAnofeventsintimestepnon}}+\VAvarover{\VAnofeventsintimestepnon}{\VAexpec{\left.\sum_{\VAeventnointimestep=0}^\VAnofeventsintimestepnon \frac{\VAdiscreteveleend_\VAeventnointimestep}{\VAscaleparam}\VAeventtimedif_\VAeventnointimestep\right|\VAnofeventsintimestepnon}}\,.\label{eq:ap_randomv0_var_lawtotalvarapplied}
\end{equation}
The same notational convention is used for the variance as for the expected value introduced in the previous section: the subscript indicates the random variable with respect to which the variance is taken and the absence of a subscript indicating that the variance is taken over all random variables, except for the ones for which it is conditioned. The second term on the right hand side of Equation~\eqref{eq:ap_randomv0_var_lawtotalvarapplied} equals zero, since the argument for $\text{Var}_{\VAnofeventsintimestepnon}$ equals the left hand side of Equation~\eqref{eq:ap_randomv0_mean_vdistriworkedout}, and was shown to be independent of $\VAnofeventsintimestepnon$ in Section~\ref{subsubsec:ap_randomv0_mean}, see Equation~\eqref{eq:ap_randomv0_mean_result}. The first term on the right hand side of Equation~\eqref{eq:ap_randomv0_var_lawtotalvarapplied} can be rewritten as a series by using Equation~\eqref{eq:ap_poissondistri}, resulting in
\begin{equation}
\VAvar{\Delta \VAdiscreteposeend}=\sum_{\VAnofeventsintimestepnon=0}^\infty\frac{1}{\VAnofeventsintimestepnon!}\left(\frac{\VAratespatial}{\VAscaleparam^2}\VAtimestep\right)^{\VAnofeventsintimestepnon}e^{-\frac{\VAratespatial}{\VAscaleparam^2}\VAtimestep}\VAvar{\left.\sum_{\VAeventnointimestep=0}^\VAnofeventsintimestepnon\frac{\VAdiscreteveleend_\VAeventnointimestep}{\VAscaleparam}\VAeventtimedif_\VAeventnointimestep\right|\VAnofeventsintimestepnon}\label{eq:ap_randomv0_var_poissonseriesapplied}
\end{equation}
We will now rewrite the conditional variance in Equation~\eqref{eq:ap_randomv0_var_poissonseriesapplied}. To do so, we first show all the time intervals $\VAeventtimedif_\VAeventnointimestep$, $\VAeventnointimestep=0\dots\VAnofeventsintimestepnon$ are identically distributed for a fixed $\VAnofeventsintimestepnon$. Then, we will be able to find a convenient formulation for the conditional variance which we will use to find an expression for the variance at the end of this section.

\paragraph{Identical distribution of the flight times in a time step.} The probability of the $\VAeventnointimestep$-th inter-event time to begin after a time in the interval $[\theta,\theta+\text{d}\theta]$ and to have a size in the interval $[\VAeventtimedif_\VAeventnointimestep,\VAeventtimedif_\VAeventnointimestep+\text{d}\VAeventtimedif_\VAeventnointimestep]$ is equal to the probability of $\VAeventnointimestep-1$ events occurring in time $\theta$ of which one occurs in the interval $[\theta,\theta+\text{d}\theta]$, multiplied by the probability of a collision occurring after exactly $\VAeventtimedif_\VAeventnointimestep$, multiplied by the probability of $\VAnofeventsintimestepnon-\VAeventnointimestep-1$ events occurring in the remaining time $\VAtimestep-\theta-\VAeventtimedif_\VAeventnointimestep$. Integrating over $\theta$, we can write the probability distribution of the time length $\VAeventtimedif_\VAeventnointimestep$ up to a constant factor as
{
\newcommand{\VAzelflocalmaxheight}{\vphantom{\frac{\left(\VArate(\VAtimestep-\VAeventtimedif_\VAeventnointimestep-\theta)\right)^{\VAnofeventsintimestepnon-\VAeventnointimestep-1}}{(\VAnofeventsintimestepnon-\VAeventnointimestep-1)!}}}
\begin{equation}
\underbrace{\VAzelflocalmaxheight e^{\!-\frac{\VAratespatial\VAeventtimedif_\VAeventnointimestep}{\VAscaleparam^2}}\!\frac{\VAratespatial}{\VAscaleparam^2}\text{d}\VAeventtimedif_\VAeventnointimestep\!\!}_{\substack{\text{one event}\\ \text{after $\VAeventtimedif_\VAeventnointimestep$}}}\int_0^{\VAtimestep-\VAeventtimedif_\VAeventnointimestep}\!\!\underbrace{\VAzelflocalmaxheight\frac{\left(\frac{\VAratespatial\theta}{\VAscaleparam^2}\right)^{\VAeventnointimestep-1}}{(\VAeventnointimestep-1)!}e^{\!-\frac{\VAratespatial\theta}{\VAscaleparam^2}}\!\frac{\VAratespatial}{\VAscaleparam^2}\text{d}\theta\!}_{\substack{\text{$\VAeventnointimestep-1$ events in $\theta$,} \\ \text{and one at $\theta$}}}\underbrace{\frac{\!\left(\frac{\VAratespatial(\VAtimestep-\VAeventtimedif_\VAeventnointimestep-\theta)}{\VAscaleparam^2}\right)^{\VAnofeventsintimestepnon-\VAeventnointimestep-1\!\!\!}}{(\VAnofeventsintimestepnon-\VAeventnointimestep-1)!}e^{\!-\frac{\VAratespatial(\VAtimestep-\VAeventtimedif_\VAeventnointimestep-\theta)}{\VAscaleparam^2}}}_\text{the remaining $\VAnofeventsintimestepnon-\VAeventnointimestep-1$ events}\,.
\end{equation}}
Restructuring this equation gives
\begin{equation}
\left(\frac{\VAratespatial}{\VAscaleparam^2}\right)^{\VAnofeventsintimestepnon} e^{-\frac{\VAratespatial\VAtimestep}{\VAscaleparam^2}}\int_0^{\VAtimestep-\VAeventtimedif_\VAeventnointimestep}\frac{\theta^{\VAeventnointimestep-1}}{(\VAeventnointimestep-1)!}\frac{(\VAtimestep-\VAeventtimedif_\VAeventnointimestep-\theta)^{\VAnofeventsintimestepnon-\VAeventnointimestep-1}}{(\VAnofeventsintimestepnon-\VAeventnointimestep-1)!}\text{d}\theta\text{d}\VAeventtimedif_\VAeventnointimestep\,,
\end{equation}
and after repeated partial integration we find
\begin{equation}
\left(\frac{\VAratespatial}{\VAscaleparam^2}\right)^{\VAnofeventsintimestepnon} e^{-\frac{\VAratespatial\VAtimestep}{\VAscaleparam^2}}\frac{(\VAtimestep-\VAeventtimedif_\VAeventnointimestep)^{\VAnofeventsintimestepnon-1}}{(\VAnofeventsintimestepnon-1)!}\text{d}\VAeventtimedif_\VAeventnointimestep\,,\label{eq:ap_randomv0_var_taudistri_intermediate}
\end{equation}
which is independent of $\VAeventnointimestep$, proving the identical distribution for each $\VAeventtimedif_\VAeventnointimestep,$ $\VAeventnointimestep\in\{0,\dots,\VAnofeventsintimestepnon\}$.

The quantity in Equation~\eqref{eq:ap_randomv0_var_taudistri_intermediate} is the probability density function of $\VAeventtimedif_\VAeventnointimestep$ up to a scaling, but it can be easily rescaled to obtain the probability density function that holds for all $\VAeventnointimestep\leq\VAnofeventsintimestepnon$:
\begin{equation}
\text{P}(\VAeventtimedif_\VAeventnointimestep|\VAnofeventsintimestepnon)=\left\{\begin{array}{ll}
\VAnofeventsintimestepnon\dfrac{(\VAtimestep-\VAeventtimedif_\VAeventnointimestep)^{\VAnofeventsintimestepnon-1}}{\VAtimestep^{\VAnofeventsintimestepnon}}&\text{if }\VAeventtimedif_\VAeventnointimestep\in[0,\VAtimestep]\\
0&\text{else}\,.
\end{array}\right.
\label{eq:ap_randomv0_var_Dtprobdens}
\end{equation}
The resulting Equation~\eqref{eq:ap_randomv0_var_Dtprobdens} holds for all $\VAeventnointimestep\in\{0,1,\dots,\VAnofeventsintimestepnon\}$.

\paragraph{Expression for the conditional variance.} To find an expression for the conditional variance in Equation~\eqref{eq:ap_randomv0_var_poissonseriesapplied}, the mean and variance of $\VAeventtimedif_\VAeventnointimestep$ will be required. These can be computed from Equation~\eqref{eq:ap_randomv0_var_Dtprobdens} to be
\begin{align}
\VAexpec{\left.\VAeventtimedif_\VAeventnointimestep\right|\VAnofeventsintimestepnon}&=\int_0^{\VAtimestep} \VAeventtimedif_\VAeventnointimestep\VAnofeventsintimestepnon\frac{(\VAtimestep-\VAeventtimedif_\VAeventnointimestep)^{\VAnofeventsintimestepnon-1}}{\VAtimestep^{\VAnofeventsintimestepnon}}\text{d}\VAeventtimedif_\VAeventnointimestep=\frac{\VAtimestep}{\VAnofeventsintimestepnon+1}\label{eq:ap_condtimestepmean}\\
\VAexpec{\left.\VAeventtimedif_\VAeventnointimestep^2\right|\VAnofeventsintimestepnon}&=\int_0^{\VAtimestep}\VAeventtimedif_\VAeventnointimestep^2\VAnofeventsintimestepnon\frac{(\VAtimestep-\VAeventtimedif_\VAeventnointimestep)^{\VAnofeventsintimestepnon-1}}{\VAtimestep^\VAnofeventsintimestepnon}\text{d}\VAeventtimedif_\VAeventnointimestep=\frac{2\VAtimestep^2}{(\VAnofeventsintimestepnon+1)(\VAnofeventsintimestepnon+2)}\label{eq:ap_condtimestepmom2}\\
\VAvar{\left.\VAeventtimedif_\VAeventnointimestep\right|\VAnofeventsintimestepnon}&=\VAexpec{\left.\VAeventtimedif_\VAeventnointimestep^2\right|\VAnofeventsintimestepnon}-\left(\VAexpec{\left.\VAeventtimedif_\VAeventnointimestep\right|\VAnofeventsintimestepnon}\right)^2=\frac{\VAnofeventsintimestepnon\VAtimestep^2}{(\VAnofeventsintimestepnon+1)^2(\VAnofeventsintimestepnon+2)}\,.\label{eq:ap_condtimestepvar}
\end{align}

Due to the fact that the different $\VAeventtimedif_\VAeventnointimestep$ and the different velocities are identically distributed, the variance of the sum formula can be applied to the conditional variance in Equation~\eqref{eq:ap_randomv0_var_poissonseriesapplied}, yielding
\begin{equation}
\VAvar{\left.\sum_{\VAeventnointimestep=0}^\VAnofeventsintimestepnon\frac{\VAdiscreteveleend_\VAeventnointimestep}{\VAscaleparam}\VAeventtimedif_\VAeventnointimestep\right|\VAnofeventsintimestepnon}=(\VAnofeventsintimestepnon+1)\left(\VAvar{\left.\frac{\VAdiscreteveleend_\VAeventnointimestep}{\VAscaleparam}\VAeventtimedif_\VAeventnointimestep\right|\VAnofeventsintimestepnon}
+\VAnofeventsintimestepnon\VAcov{\left.\frac{\VAdiscreteveleend_\VAeventnointimestep}{\VAscaleparam}\VAeventtimedif_\VAeventnointimestep\right|\VAnofeventsintimestepnon}\right)\label{eq:ap_randomv0_var_conddxvar_sumformula}
\end{equation}

Because of the independence of the $\VAdiscreteveleend_\VAeventnointimestep$ to other random variables, we have
\begin{multline}
\VAvar{\left.\frac{\VAdiscreteveleend_\VAeventnointimestep}{\VAscaleparam}\VAeventtimedif_\VAeventnointimestep\right|\VAnofeventsintimestepnon}=\VAvar{\left.\frac{\VAdiscreteveleend_\VAeventnointimestep}{\VAscaleparam}\right|\VAnofeventsintimestepnon}\VAvar{\left.\VAeventtimedif_\VAeventnointimestep\right|\VAnofeventsintimestepnon}
+\VAvar{\left.\frac{\VAdiscreteveleend_\VAeventnointimestep}{\VAscaleparam}\right|\VAnofeventsintimestepnon}\VAexpec{\left.\VAeventtimedif_\VAeventnointimestep\right|\VAnofeventsintimestepnon}^2\\
+\VAvar{\left.\VAeventtimedif_\VAeventnointimestep\right|\VAnofeventsintimestepnon}\VAexpec{\left.\frac{\VAdiscreteveleend_\VAeventnointimestep}{\VAscaleparam}\right|\VAnofeventsintimestepnon}^2
\label{eq:ap_randomv0_var_condvdtvar_splitup}
\end{multline}
and
\begin{align}
\VAcov{\left.\frac{\VAdiscreteveleend_\VAeventnointimestep}{\VAscaleparam}\VAeventtimedif_\VAeventnointimestep\right|\VAnofeventsintimestepnon}&=\VAcov{\VAeventtimedif_\VAeventnointimestep|\VAnofeventsintimestepnon}\VAexpec{\left.\frac{\VAdiscreteveleend_\VAeventnointimestep}{\VAscaleparam}\right|\VAnofeventsintimestepnon}^2\,.\label{eq:ap_randomv0_var_condvtcov}
\end{align}
Applying the standard formula for the variance of a sum formula to the event times, we find
\begin{equation}
\VAvar{\left.\sum_{\VAeventnointimestep=0}^\VAnofeventsintimestepnon\VAeventtimedif_\VAeventnointimestep\right|\VAnofeventsintimestepnon}=(\VAnofeventsintimestepnon+1)\VAvar{\VAeventtimedif_\VAeventnointimestep|\VAnofeventsintimestepnon}+(\VAnofeventsintimestepnon+1)\VAnofeventsintimestepnon\VAcov{\VAeventtimedif_\VAeventnointimestep|\VAnofeventsintimestepnon}=0\,,
\end{equation}
which allows the expression of the covariance of the flight time contributions from Equation~\eqref{eq:ap_randomv0_var_condvtcov} as
\begin{equation}
\VAcov{\left.\VAeventtimedif_\VAeventnointimestep\frac{\VAdiscreteveleend_\VAeventnointimestep}{\VAscaleparam}\right|\VAnofeventsintimestepnon}=-\frac{\VAvar{\VAeventtimedif_\VAeventnointimestep|\VAnofeventsintimestepnon}}{\VAnofeventsintimestepnon}\VAexpec{\left.\frac{\VAdiscreteveleend_\VAeventnointimestep}{\VAscaleparam}\right|\VAnofeventsintimestepnon}^2\,.\label{eq:ap_randomv0_var_condvdtcov}
\end{equation}
Equations~\eqref{eq:ap_randomv0_var_condvdtvar_splitup} and~\eqref{eq:ap_randomv0_var_condvdtcov} can be used to transform the conditional variance of Equation~\eqref{eq:ap_randomv0_var_poissonseriesapplied} to
\begin{equation}
\VAvar{\left.\sum_{\VAeventnointimestep=0}^\VAnofeventsintimestepnon\frac{\VAdiscreteveleend_\VAeventnointimestep}{\VAscaleparam}\VAeventtimedif_\VAeventnointimestep\right|\VAnofeventsintimestepnon}=(\VAnofeventsintimestepnon+1)\VAvar{\left.\frac{\VAdiscreteveleend_\VAeventnointimestep}{\VAscaleparam}\right|\VAnofeventsintimestepnon}\left(\VAvar{\left.\VAeventtimedif_\VAeventnointimestep\right|\VAnofeventsintimestepnon}
+\VAexpec{\left.\VAeventtimedif\right|\VAnofeventsintimestepnon}^2\right)\,.\label{eq:ap_randomv0_var_conddxvar_intermediate}
\end{equation}
All the terms of Equation~\eqref{eq:ap_randomv0_var_conddxvar_intermediate} are known: the mean and variance of $\VAeventtimedif_\VAeventnointimestep$ conditioned on $\VAnofeventsintimestepnon$ are given by Equations~\eqref{eq:ap_condtimestepmean} and~\eqref{eq:ap_condtimestepvar}, and the velocities are independent of $\VAnofeventsintimestepnon$ and have variance $\VAbackgroundspeedtemp$. Using this knowledge, transforms Equation~\eqref{eq:ap_randomv0_var_conddxvar_intermediate} into
\begin{equation}
\VAvar{\left.\sum_{\VAeventnointimestep=0}^\VAnofeventsintimestepnon\frac{\VAdiscreteveleend_\VAeventnointimestep}{\VAscaleparam}\VAeventtimedif_\VAeventnointimestep\right|\VAnofeventsintimestepnon}=\frac{\VAbackgroundspeedtemp\VAnofeventsintimestepnon\VAtimestep^2}{\VAscaleparam^2(\VAnofeventsintimestepnon+1)(\VAnofeventsintimestepnon+2)}
+\frac{\VAbackgroundspeedtemp\VAtimestep^2}{\VAscaleparam^2(\VAnofeventsintimestepnon+1)}=\frac{\VAbackgroundspeedtemp}{\VAscaleparam^2}\VAtimestep^2\frac{2}{\VAnofeventsintimestepnon+2\,.}\label{eq:ap_randomv0_conddxvar}
\end{equation}

\paragraph{Expression for the variance.} Expression~\eqref{eq:ap_randomv0_conddxvar} for the conditional variance of Equation~\eqref{eq:ap_randomv0_var_poissonseriesapplied} can be used to obtain
\begin{align}
\VAvar{\Delta\VAdiscreteposeend}&=\sum_{\VAnofeventsintimestepnon=0}^\infty \frac{\VAbackgroundspeedtemp}{\VAscaleparam^2}\VAtimestep^2e^{-\frac{\VAratespatial\VAtimestep}{\VAscaleparam^2}}\frac{2}{\VAnofeventsintimestepnon+2}\frac{\left(\frac{\VAratespatial}{\VAscaleparam^2}\VAtimestep\right)^\VAnofeventsintimestepnon}{\VAnofeventsintimestepnon!}\\
&=2\frac{\VAbackgroundspeedtemp}{\VAscaleparam^2}\VAtimestep^2e^{-\frac{\VAratespatial}{\VAscaleparam^2}\VAtimestep}\sum_{\VAnofeventsintimestepnon=0}^\infty(\VAnofeventsintimestepnon+1)\frac{\left(\frac{\VAratespatial\VAtimestep}{\VAscaleparam^2}\right)^{\VAnofeventsintimestepnon}}{(\VAnofeventsintimestepnon+2)!}\,,
\intertext{which can be rewritten in the form of the series expansion of an exponential and its derivative as}
\VAvar{\Delta\VAdiscreteposeend}&=2\frac{\VAbackgroundspeedtemp}{\VAscaleparam^2}\VAtimestep^2e^{-\frac{\VAratespatial\VAtimestep}{\VAscaleparam^2}}\left(\sum_{\VAnofeventsintimestepnon=0}^\infty(\VAnofeventsintimestepnon+2)\frac{\left(\frac{\VAratespatial\VAtimestep}{\VAscaleparam^2}\right)^{\VAnofeventsintimestepnon}}{(\VAnofeventsintimestepnon+2)!}-\sum_{\VAnofeventsintimestepnon=0}^\infty\frac{\left(\frac{\VAratespatial\VAtimestep}{\VAscaleparam^2}\right)^{\VAnofeventsintimestepnon}}{(\VAnofeventsintimestepnon+2)!}\right)\\
&=2\frac{\VAbackgroundspeedtemp}{\VAscaleparam^2}\VAtimestep^2e^{-\frac{\VAratespatial\VAtimestep}{\VAscaleparam^2}}\left(\frac{\VAscaleparam^2}{\VAratespatial\VAtimestep}\sum_{\VAnofeventsintimestepnon=2}^\infty\VAnofeventsintimestepnon\frac{\left(\frac{\VAratespatial\VAtimestep}{\VAscaleparam^2}\right)^{\VAnofeventsintimestepnon-1}}{\VAnofeventsintimestepnon!}-\left(\frac{\VAscaleparam^2}{\VAratespatial\VAtimestep}\right)^2\sum_{\VAnofeventsintimestepnon=2}^\infty\frac{\left(\frac{\VAratespatial\VAtimestep}{\VAscaleparam^2}\right)^{\VAnofeventsintimestepnon}}{\VAnofeventsintimestepnon}\right)\\
&=2\frac{\VAscaleparam^2}{\VAratespatial^2}\VAbackgroundspeedtemp e^{-\frac{\VAratespatial\VAtimestep}{\VAscaleparam^2}}\left(\frac{\VAratespatial\VAtimestep}{\VAscaleparam^2}\frac{\text{d}\left(e^{\frac{\VAratespatial\VAtimestep}{\VAscaleparam^2}}-1-\frac{\VAratespatial\VAtimestep}{\VAscaleparam^2}\right)}{\text{d}\frac{\VAratespatial\VAtimestep}{\VAscaleparam^2}}-\left(e^{\frac{\VAratespatial\VAtimestep}{\VAscaleparam^2}}-1-\frac{\VAratespatial\VAtimestep}{\VAscaleparam^2}\right)\right)\,,
\end{align}
resulting in the formula for the variance
\begin{equation}
\VAvar{\Delta\VAdiscreteposeend}=2\frac{\VAscaleparam^2}{\VAratespatial^2}\VAbackgroundspeedtemp\left(e^{-\frac{\VAratespatial\VAtimestep}{\VAscaleparam^2}}-1+\frac{\VAratespatial\VAtimestep}{\VAscaleparam^2}\right)\,.\label{eq:ap_randomv0_var}
\end{equation}
In the diffusive limit, $\VAscaleparam\rightarrow0$, the variance in Equation~\eqref{eq:ap_randomv0_var} indeed becomes equal to $2\frac{\VAbackgroundspeedtemp}{\VAratespatial}\VAtimestep$, in correspondence with Equation~\eqref{eq:ap_modelsim_difflimit_sdeequiv} of Section~\ref{subsec:ap_modelsim_difflimit}. On the other hand, for fixed $\VAscaleparam$ and decreasing time step size, we obtain a variance proportional to $\frac{\VAbackgroundspeedtemp}{\VAscaleparam^2}\VAtimestep^2$. This equals the variance of the motion during a time $\VAtimestep$ where the particle has a constant but randomly sampled velocity with variance $\frac{\VAbackgroundspeedtemp}{\VAscaleparam^2}$, which is exactly the kinetic limit of the particle process, with velocities $\frac{v}{\VAscaleparam}$ with $v$ distributed according to the Maxwellian of Equation~\eqref{eq:ap_kinetic_postcolveldistr_pdf}.

\subsection{Including correlation by conditioning on the final velocity\label{subsec:ap_fixedv0}}

In Section~\ref{subsec:ap_randomv0}, all velocities were treated as independent new samples. When multiple time steps are considered, the last velocity of the previous time step is equal to the first of the next time step, which results in correlation of the positional increment in different time steps. To correlate the motion within a time step to its subsequent time steps, the mean and variance of the motion can take into account the final velocity during the time step. In a low-collisional regime, this correlation is an important feature of the kinetic behaviour.

\subsubsection{Including correlation by conditioning on the final velocity: mean}

As in the previous section, we condition the expected value on the number of events and can use the independence of the velocities with respect to the other variables, resulting in
\begin{align}
\VAexpec{\Delta \VAdiscreteposeend}&=\VAexpec{\VAexpec{\left.\sum_{\VAeventnointimestep=0}^\VAnofeventsintimestepnon\VAexpec{\frac{\VAdiscreteveleend_\VAeventnointimestep}{\VAscaleparam}}\VAeventtimedif_\VAeventnointimestep\right|\VAnofeventsintimestepnon}}\,.
\end{align}
The difference with before is the treatment of $\VAdiscreteveleend_\VAnofeventsintimestepnon$ as a fixed predetermined value, resulting in
\begin{align}
\VAexpec{\Delta \VAdiscreteposeend}&=\VAexpec{\VAexpec{\left.\sum_{\VAeventnointimestep=0}^{\VAnofeventsintimestepnon-1}\VAbackgroundspeedmeanflow\VAeventtimedif_\VAeventnointimestep+\frac{\VAdiscreteveleend_\VAnofeventsintimestepnon}{\VAscaleparam}\VAeventtimedif_\VAnofeventsintimestepnon\right|\VAnofeventsintimestepnon}}\\
&=\VAexpec{\VAexpec{\left.\VAbackgroundspeedmeanflow(\VAtimestep-\VAeventtimedif_\VAnofeventsintimestepnon)+\frac{\VAdiscreteveleend_\VAnofeventsintimestepnon}{\VAscaleparam}\VAeventtimedif_\VAnofeventsintimestepnon\right|\VAnofeventsintimestepnon}}\\
&=\VAexpec{\VAbackgroundspeedmeanflow(\VAtimestep-\VAeventtimedif_\VAnofeventsintimestepnon)+\frac{\VAdiscreteveleend_\VAnofeventsintimestepnon}{\VAscaleparam}\VAeventtimedif_\VAnofeventsintimestepnon}\\
&=\VAbackgroundspeedmeanflow\VAtimestep+\left(\frac{\VAdiscreteveleend_\VAnofeventsintimestepnon}{\VAscaleparam}-\VAbackgroundspeedmeanflow\right)\VAexpec{\VAeventtimedif_\VAnofeventsintimestepnon}\,.
\end{align}
Since $\VAeventtimedif_\VAnofeventsintimestepnon$ equals the time before any collision occurs in $\VAtimestep$, it is exponentially distributed with rate $\frac{\VAratespatial}{\VAscaleparam^2}$, and with $\VAtimestep$ as a ceiling. Its expected value can consequently be computed as
\begin{align}
\VAexpec{\VAeventtimedif_\VAnofeventsintimestepnon}&=\int_0^{\VAtimestep}\VAeventtimedif\frac{\VAratespatial}{\VAscaleparam^2} e^{-\frac{\VAratespatial\VAeventtimedif}{\VAscaleparam^2}}\text{d}\VAeventtimedif+\VAtimestep\int_{\VAtimestep}^\infty\frac{\VAratespatial}{\VAscaleparam^2} e^{-\frac{\VAratespatial}{\VAscaleparam^2}\VAeventtimedif}\text{d}\VAeventtimedif\\
&=\frac{\VAscaleparam^2}{\VAratespatial}-\frac{\VAscaleparam^2}{\VAratespatial} e^{-\frac{\VAratespatial\VAtimestep}{\VAscaleparam^2}}-\VAtimestep e^{-\frac{\VAratespatial\VAtimestep}{\VAscaleparam^2}}+\VAtimestep e^{-\frac{\VAratespatial\VAtimestep}{\VAscaleparam^2}}\\
&=\frac{\VAscaleparam^2}{\VAratespatial}\left(1-e^{-\frac{\VAratespatial\VAtimestep}{\VAscaleparam^2}}\right)\,,\label{eq:ap_fixedv0_dt0mean}
\end{align}
resulting in the mean
\begin{equation}
\VAexpec{\Delta\VAdiscreteposeend}=\VAbackgroundspeedmeanflow\VAtimestep+\left(\frac{\VAdiscreteveleend_\VAnofeventsintimestepnon}{\VAscaleparam}-\VAbackgroundspeedmeanflow\right)\frac{\VAscaleparam^2}{\VAratespatial}\left(1-e^{-\frac{\VAratespatial\VAtimestep}{\VAscaleparam^2}}\right)\,.\label{eq:ap_fixedv0_mean_result}
\end{equation}
In the diffusive limit, the mean becomes equal to $\VAbackgroundspeedmeanflow\VAtimestep$. Conversely, for fixed $\VAscaleparam$ and $\VAtimestep\rightarrow0$, the mean becomes equal to $\VAdiscreteveleend_{\VAnofeventsintimestepnon}$ since the probability of a velocity changes becomes zero.

\subsubsection{Including correlation by conditioning on the final velocity: variance}

As in Section~\ref{subsubsec:ap_fixedv0_var}, we start with the law of total variance, but now applied to the last time interval $\VAeventtimedif_\VAnofeventsintimestepnon$, giving
\begin{equation}
\VAvar{\Delta\VAdiscreteposeend}=\VAexpec{\VAvar{\left.\sum_{\VAeventnointimestep=0}^\VAnofeventsintimestepnon\frac{\VAdiscreteveleend_\VAeventnointimestep}{\VAscaleparam}\VAeventtimedif_\VAeventnointimestep\right|\VAeventtimedif_\VAnofeventsintimestepnon}}+\VAvar{\VAexpec{\left.\sum_{\VAeventnointimestep=0}^\VAnofeventsintimestepnon\frac{\VAdiscreteveleend_\VAeventnointimestep}{\VAscaleparam}\VAeventtimedif_\VAeventnointimestep\right|\VAeventtimedif_\VAnofeventsintimestepnon}}\,.
\label{eq:ap_fixedv0_var_lawtotalvarapplied}
\end{equation}
For the first term of Equation~\eqref{eq:ap_fixedv0_var_lawtotalvarapplied}, we can use that $V_\VAnofeventsintimestepnon$ is a known and fixed value to write
\begin{align}
\VAexpec{\VAvar{\left.\sum_{\VAeventnointimestep=0}^\VAnofeventsintimestepnon\frac{\VAdiscreteveleend_\VAeventnointimestep}{\VAscaleparam}\VAeventtimedif_\VAeventnointimestep\right|\VAeventtimedif_\VAnofeventsintimestepnon}}&=\VAexpec{\VAvar{\left.\sum_{\VAeventnointimestep=0}^{\VAnofeventsintimestepnon-1}\frac{\VAdiscreteveleend_\VAeventnointimestep}{\VAscaleparam}\VAeventtimedif_\VAeventnointimestep\right|\VAeventtimedif_\VAnofeventsintimestepnon}}\,.\label{eq:ap_fixedv0_var_lawtotalvar_term1}
\end{align}
Equation~\eqref{eq:ap_fixedv0_var_lawtotalvar_term1} conforms to the case in Section~\ref{subsec:ap_randomv0} but with a time step of length $\VAtimestep-\VAeventtimedif_\VAnofeventsintimestepnon$, enabling the use of Equation~\eqref{eq:ap_randomv0_var} to find
\begin{equation}
\VAexpec{\VAvar{\left.\sum_{\VAeventnointimestep=0}^\VAnofeventsintimestepnon\frac{\VAdiscreteveleend_\VAeventnointimestep}{\VAscaleparam}\VAeventtimedif_\VAeventnointimestep\right|\VAeventtimedif_\VAnofeventsintimestepnon}}=\VAexpec{2\frac{\VAscaleparam^2}{\VAratespatial^2}\VAbackgroundspeedtemp\left(e^{-\frac{\VAratespatial\left(\VAtimestep-\VAeventtimedif_\VAnofeventsintimestepnon\right)}{\VAscaleparam^2}}-1+\frac{\VAratespatial(\VAtimestep-\VAeventtimedif_\VAnofeventsintimestepnon)}{\VAscaleparam^2}\right)}\,,
\end{equation}
\begin{align}
\hphantom{jajaja}&=\int_0^{\VAtimestep}2\frac{\VAscaleparam^2}{\VAratespatial^2}\VAbackgroundspeedtemp\left(e^{-\frac{\VAratespatial(\VAtimestep-\VAeventtimedif_\VAnofeventsintimestepnon)}{\VAscaleparam^2}}-1+\frac{\VAratespatial(\VAtimestep-\VAeventtimedif_\VAnofeventsintimestepnon)}{\VAscaleparam^2}\right)\frac{\VAratespatial}{\VAscaleparam^2} e^{-\frac{\VAratespatial(\VAeventtimedif_\VAnofeventsintimestepnon)}{\VAscaleparam}}\text{d}\VAeventtimedif_\VAnofeventsintimestepnon\\
&\ \ \ \ \ \ \ +\left.2\frac{\VAscaleparam^2}{\VAratespatial^2}\VAbackgroundspeedtemp\left(e^{-\frac{\VAratespatial(\VAtimestep-\VAeventtimedif_\VAnofeventsintimestepnon)}{\VAscaleparam^2}}-1+\frac{\VAratespatial(\VAtimestep-\VAeventtimedif_\VAnofeventsintimestepnon)}{\VAscaleparam^2}\right)\right|_{\VAeventtimedif_\VAnofeventsintimestepnon=\VAtimestep}e^{-\frac{\VAratespatial\VAtimestep}{\VAscaleparam^2}}\\
&=2\frac{\VAscaleparam^2}{\VAratespatial^2}\VAbackgroundspeedtemp\left(2e^{-\frac{\VAratespatial\VAtimestep}{\VAscaleparam^2}}+\frac{\VAratespatial\VAtimestep}{\VAscaleparam^2}+\frac{\VAratespatial\VAtimestep}{\VAscaleparam^2} e^{-\frac{\VAratespatial\VAtimestep}{\VAscaleparam^2}}-2\right)\,.\label{eq:ap_fixedv0_lawoftotalvar_term1}
\end{align}
The second term in Equation~\eqref{eq:ap_fixedv0_var_lawtotalvarapplied} can be reworked as
\begin{align}
\VAvar{\VAexpec{\left.\sum_{\VAeventnointimestep=0}^\VAnofeventsintimestepnon\frac{\VAdiscreteveleend_\VAeventnointimestep}{\VAscaleparam}\VAeventtimedif_\VAeventnointimestep\right|\VAeventtimedif_\VAnofeventsintimestepnon}}&=\VAvar{\frac{\VAdiscreteveleend_\VAnofeventsintimestepnon}{\VAscaleparam}\VAeventtimedif_\VAnofeventsintimestepnon+\VAbackgroundspeedmeanflow(\VAtimestep-\VAeventtimedif_\VAnofeventsintimestepnon)}\\
&=\left(\frac{\VAdiscreteveleend_\VAnofeventsintimestepnon}{\VAscaleparam}-\VAbackgroundspeedmeanflow\right)^2\VAvar{\VAeventtimedif_\VAnofeventsintimestepnon}\,,\label{eq:ap_fixedv0_lawoftotalvar_term2}
\end{align}
by using the independence of the $\VAdiscreteveleend_\VAeventnointimestep$ and the fact that the time intervals sum to $\VAtimestep$. For the variance on $\VAeventtimedif_{\VAnofeventsintimestepnon}$, we use Equation~\eqref{eq:ap_fixedv0_dt0mean} and partial integration to find
\begin{align}
\VAvar{\VAeventtimedif_{\VAnofeventsintimestepnon}}&=\int_0^{\VAtimestep}\frac{\VAratespatial\VAeventtimedif_\VAnofeventsintimestepnon^2}{\VAscaleparam^2} e^{-\frac{\VAratespatial\VAeventtimedif_\VAnofeventsintimestepnon}{\VAscaleparam^2}}\text{d}\VAeventtimedif_\VAnofeventsintimestepnon+\int_{\VAtimestep}^{\infty}\frac{\VAratespatial\VAtimestep^2}{\VAscaleparam^2} e^{-\frac{\VAratespatial\VAeventtimedif_\VAnofeventsintimestepnon}{\VAscaleparam^2}}\text{d}\VAeventtimedif_\VAnofeventsintimestepnon-\left(\frac{\VAscaleparam^2}{\VAratespatial}\left(1-e^{-\frac{\VAratespatial\VAtimestep}{\VAscaleparam^2}}\right)\right)^2\\
&=\frac{\VAscaleparam^4}{\VAratespatial^2}\left(1-2\frac{\VAratespatial\VAtimestep}{\VAscaleparam^2} e^{-\frac{\VAratespatial\VAtimestep}{\VAscaleparam^2}}-e^{-2\frac{\VAratespatial\VAtimestep}{\VAscaleparam^2}}\right)\,.\label{eq:ap_fixedv0_vardt0}
\end{align}
Using the result of Equation~\eqref{eq:ap_fixedv0_vardt0} in Equation~\eqref{eq:ap_fixedv0_lawoftotalvar_term2} yields an expression the second term of Equation~\eqref{eq:ap_fixedv0_var_lawtotalvarapplied}. Together with the first term, as expressed by Equation~\eqref{eq:ap_fixedv0_lawoftotalvar_term1}, we find the expression for the variance of a positional increment conditioned on the final velocity:
\begin{multline}
\VAvar{\Delta \VAdiscreteposeend}=2\VAbackgroundspeedtemp\frac{\VAscaleparam^2}{\VAratespatial^2}\left(2e^{-\frac{\VAratespatial\VAtimestep}{\VAscaleparam^2}}+\frac{\VAratespatial\VAtimestep}{\VAscaleparam^2}+\frac{\VAratespatial\VAtimestep}{\VAscaleparam^2} e^{-\frac{\VAratespatial\VAtimestep}{\VAscaleparam^2}}-2\right)\\
+\left(\frac{\VAdiscreteveleend_\VAnofeventsintimestepnon}{\VAscaleparam}-\VAbackgroundspeedmeanflow\right)^2\frac{\VAscaleparam^4}{\VAratespatial^2}\left(1-2\frac{\VAratespatial}{\VAscaleparam^2}\VAtimestep e^{-\frac{\VAratespatial\VAtimestep}{\VAscaleparam^2}}-e^{-2\frac{\VAratespatial\VAtimestep}{\VAscaleparam^2}}\right)\,.\label{eq:ap_fixedv0_var}
\end{multline}
In the diffusive limit, $\VAscaleparam\rightarrow0$, Equation~\eqref{eq:ap_fixedv0_var} becomes equal to $2\frac{\VAbackgroundspeedtemp}{\VAratespatial}\VAtimestep$. This is in accordance with Equation~\eqref{eq:ap_modelsim_difflimit_sderandomwalk} and is identical to the diffusive limit for Equation~\eqref{eq:ap_randomv0_var}. For fixed $\VAscaleparam$ and decreasing $\VAtimestep$, Equation~\eqref{eq:ap_fixedv0_var} now becomes
\begin{equation}
\left(\VAbackgroundspeedtemp+\left(V_K-\VAscaleparam u\right)^2\right)\frac{\VAratespatial}{3\VAscaleparam^4}\VAtimestep^3\,,\ \ \VAtimestep\rightarrow0\,.
\end{equation}
This limiting variance is of higher order in $\VAtimestep$ when $\VAtimestep\rightarrow0$, than for Equation~\eqref{eq:ap_randomv0_var}, where the final velocity was not fixed. This means that with fixed final velocity, the variance becomes lower, which is to be expected since part of the randomness has been removed by fixing the final velocity.

\section{Kinetic-diffusion scheme\label{sec:ap_newscheme}}

In this Section, we introduce the modified simulation strategy that has been announced in Section~\ref{subsec:ap_modelsim_strategy} in which the many free flights in the high-collisional regime are replaced by a single random walk step of the form given by Equation~\eqref{eq:ap_modelsim_strategy_randomwalk}. This replacement by a random walk step is expected to function well when the collision rate is very high, since an advection-diffusion process is the limit of the kinetic process in the high-collisional limit as shown in Section~\ref{subsec:ap_modelsim_difflimit}. In lower-collisional regimes on the other hand, an advection-diffusion description of the particle motion is invalid. We will design the algorithm in such a way that it also functions well in a low-collisional setting.

A first important error when moving diffusively with a random walk step instead of kinetically is that the next position will falsely become normally distributed. A second important error is the supposed independence of the motion during subsequent time steps by a random walk. When moving according to the kinetic model, subsequent time steps are correlated due to the final velocity in a time step featuring as the initial velocity in the next time step. To maintain the kinetic nature of the motion when the collisionality is low, we combine two different strategies in our algorithm.

The first strategy consists of executing the first flight during a time step kinetically with as the initial velocity, correctly, the final velocity of the previous time step. Only if a collision occurs within the time step, the remainder, $\VAdiftime\leq\VAtimestep$, is filled with a random walk step. This way, if the probability of a collision occurring in a time step is low, the particle will nearly always move according to the kinetic process, and only for a small fraction of the time steps, partly diffusively.

The second strategy improves the situation further by incorporating the correlation with the next time step in the advection and diffusion coefficient of the random walk movement for a time $\VAdiftime$. To do so, the background parameters, $\VAdiftime$, $\VAratespatial$, $\VAbackgroundspeedmeanflow$, and $\VAbackgroundspeedtemp$ are assumed to be spatially homogeneous. Then, the results from Section~\ref{subsec:ap_fixedv0} apply and the exact mean and variance of the kinetic process conditioned on the final velocity equalling the initial velocity of the next time step can be used. We mimic Equation~\eqref{eq:ap_modelsim_strategy_randomwalk}, but use the advection term
\begin{equation}
A_\VAscaleparam\VAdiftime=\VAexpec{\Delta \VAdiscreteposeend}\,,
\end{equation}
with $\VAexpec{\Delta \VAdiscreteposeend}$ from Equation~\eqref{eq:ap_fixedv0_mean_result} and the diffusion term
\begin{equation}
\sqrt{2D_\VAscaleparam\VAdiftime}=\sqrt{\VAvar{\Delta \VAdiscreteposeend}}\,,
\end{equation}
with $\VAvar{\Delta \VAdiscreteposeend}$ as in Equation~\eqref{eq:ap_fixedv0_var}. For the time step in these equations, we use the time $\VAdiftime$ during which the particle moves diffusively. For the other parameters $\VAratespatial$, $\VAbackgroundspeedmeanflow$, and $\VAbackgroundspeedtemp$ that appear in Equations~\eqref{eq:ap_fixedv0_mean_result} and~\eqref{eq:ap_fixedv0_var}, we use the local values at the beginning of the diffusive motion. For the final velocity in the time step, we use a fresh sample, which is then used as the initial velocity in the next time step. To cope with highly heterogeneous $\VAratespatial$, some modifications are presented in~\cite{mortier2019KDfusioncase}.

The resulting algorithm is shown as Algorithm~\ref{alg:ap_kindifeNe}, which we refer to as KD.

\newcommand\mycommfont[1]{\footnotesize\ttfamily\textcolor{blue}{#1}}
\begin{algorithm}[H]
  \caption{A kinetic-diffusion simulation up to time $\VAtotaltime=\VAnoftimesteps\VAtimestep$\label{alg:ap_kindifeNe}}
  \begin{algorithmic}[1]
    \Function{DiffusionKinetic\_KD}{$\VAdiscreteposeend,\VAdiscreteveleend,\VAtimestep,\VAratespatial(x),\VAscaleparam,\VApostcolveldistri(v;x),\VAbackgroundspeedmeanflow(x),\VAbackgroundspeedtemp(x)$}
    	\State $t\leftarrow0$
    	\While{$\VAtimevar<\VAtotaltime$}
		    \State $\VAeventtimedif\gets$ \Call{SampleCollision}{$\VAdiscreteposeend,\VAdiscreteveleend,\VAratespatial(\VAdiscreteposeend),\VAscaleparam$}\Comment{see page~\pageref{plaats:ap_kin_eventtime_function}}
		    \State $\tau\gets$\Call{min}{$\VAeventtimedif,\VAtotaltime-t$}\Comment{determine the kinetically moved time}
		    \State $X\gets X+\frac{V}{\VAscaleparam}\tau$\Comment{kinetic part of the motion}
			\If{$\VAeventtimedif<\VAtotaltime-\VAtimevar$}\Comment{check if a collision occurred}
          		\State $\VAdiscreteveleend\gets\VApostcolveldistri(v;x)$\Comment{sample the velocity for the next time step}
			    \State $\theta\gets \VAtimestep-(\tau$ \textproc{mod} $\VAtimestep)$\Comment{determine the remainder of the time step}
          		\State $\VAdiscreteposeend\gets \VAdiscreteposeend+\left.\mathcal{N}\left(\text{eq.}\eqref{eq:ap_fixedv0_mean_result},\text{eq.}\eqref{eq:ap_fixedv0_var}\right)\right|_{\begin{subarray}{l}\VAtimestep=\VAtimestep-\VAeventtimedif,\VAdiscreteveleend_\VAnofeventsintimestepnon=\VAdiscreteveleend,\VAscaleparam=\VAscaleparam\\ \VAratespatial=\VAratespatial(\VAdiscreteposeend),\VAbackgroundspeedmeanflow=\VAbackgroundspeedmeanflow(\VAdiscreteposeend),\VAbackgroundspeedtemp=\VAbackgroundspeedtemp(\VAdiscreteposeend)\end{subarray}}$\label{alg:ap_kindifeNe_D}\Comment{random walk}
	        \EndIf
		    \State $t\rightarrow t+\VAtimestep$
		\EndWhile
      \State \Return{$\VAdiscreteposeend,\VAdiscreteveleend$}
    \EndFunction
  \end{algorithmic}
\end{algorithm}
When the collisionality $\frac{\VAratespatial}{\VAscaleparam^2}\VAtimestep$ is large, Algorithm~\ref{alg:ap_kindifeNe} has an enormous computational advantage over Algorithm~\ref{alg:ap_kin} in that it does not require the execution of an expected $\frac{\VAratespatial}{\VAscaleparam^2}\VAtimestep$ collisions in a time step but at most one. In low-collisional regimes, the algorithm collapses to the standard Monte Carlo method, meaning approximately the same number of collisions occur.

In the next two sections, we will present an analysis of the error of the KD scheme. The result of Section~\ref{sec:ap_limnul} will show Algorithm~\ref{alg:ap_kindifeNe} converges to the kinetic algorithm of Section~\ref{sec:ap_kinetic} when the scaling parameter $\VAscaleparam$ is finite and the time step decreases. In Section~\ref{sec:ap_liminf}, we will show the error vanishes for $\VAscaleparam\rightarrow0$ as well. Then, in Section~\ref{sec:ap_num}, the low error and reduction in computational time are illustrated numerically.

\section{Consistency and convergence analysis for finite values of the scaling parameter\label{sec:ap_limnul}}

In this section, we analyze the error by replacing a kinetic simulation by the algorithm proposed in Section~\ref{sec:ap_newscheme} when the scaling parameter $\VAscaleparam$ is finite. In this regime, the consistency and convergence analyses require only an analysis of how the error depends on the time step $\VAtimestep$. We conduct this analysis for a spatially homogeneous setting, where $\VAratespatial(x)\equiv\VAratespatial$, $\VAbackgroundspeedmeanflow(x)\equiv\VAbackgroundspeedmeanflow$, and $\VAbackgroundspeedtemp(x)\equiv\VAbackgroundspeedtemp$. In Section~\ref{sec:ap_liminf}, we will treat the diffusive limit, when $\VAscaleparam\rightarrow0$, separately.

We evaluate the error in terms of its Wasserstein metric~\cite{villani2008optimaltransportWasserstein}. The Wasserstein distance between two distributions $g(x)$ and $h(x)$ considers the distance between two distributions as the minimal distance the mass needs to be moved on average to obtain one distribution ($g(x)$) from the other ($h(x)$). The Wasserstein distance can be expressed as the Wasserstein norm of the difference between the distributions $g(x)-h(x)$,
\begin{equation}
\zelfkalk{g(x)-h(x)}=\inf_{\mathcal{J}(g,h)}\iint|x-y|j(x,y)\text{d}x\text{d}y\,,\label{eq:ap_intro_W1def}
\end{equation}
with $\mathcal{J}(g,h)$ the class of all joint probabilities $j(x,y)$ with marginal distributions $g(x)$ and $h(y)$. An important feature of the Wasserstein distance is the consideration of the distance of mass differences, which is relevant when considering the physical meaning of particle distributions. Other often used distances like those derived from the L$_p$-norms penalize mass difference only, regardless of how far the mass of the other distribution is located.

We will compare the exact distribution at time $t^{(\VAtimestepno+1)}=(\VAtimestepno+1)\VAtimestep$, denoted as $f(t^{(n+1)})$, with the distribution $f^{(n+1)}$ obtained by the kinetic-diffusion simulation. We write $\VAaprealoperator_{\VAtimestep}$ for the operator that evolves a distribution according to the kinetic equation, Equation~\eqref{eq:ap_kinetic_integrodiff}, over a time interval of length $\VAtimestep$ and $\VAapkindifoperator_{\VAtimestep}$ for the operator that corresponds to the hybridized KD simulation scheme with time step $\VAtimestep$. With these notations, we can write
\begin{equation}
f(t^{(n+1)})=f(t^{(n)})+\VAaprealoperator_{\VAtimestep}(f(t^{(n)}))\quad\text{and}\quad f^{(n+1)}=f^{(n)}+\VAapkindifoperator_{\VAtimestep}(f^{(n)})\,.
\end{equation}
The difference between both is thus
\begin{equation}
f(t^{(n+1)})-f^{(n+1)}=f(t^{(n)})-f^{(n)}+\VAaprealoperator_{\VAtimestep}(f(t^{(n)}))-\VAapkindifoperator_{\VAtimestep}(f^{(n)})\,.\label{eq:ap_limnul_distridif_simple}
\end{equation}
To the right hand side of~\eqref{eq:ap_limnul_distridif_simple}, we add and subtract the term $\VAaprealoperator_{\VAtimestep}(f^{(\VAtimestepno)})$, and we can then bound the Wasserstein error on the distribution via the triangle inequality as
\begin{multline}
\zelfkalkdun{f(t^{(\VAtimestepno+1)})-f^{(\VAtimestepno+1)}}\leq\underbrace{\zelfkalkdun{f(t^{(\VAtimestepno)})+\VAaprealoperator_{\VAtimestep}(f(t^{(\VAtimestepno)}))-\left(f^{(\VAtimestepno)}+\VAaprealoperator_{\VAtimestep}(f^{(\VAtimestepno)})\right)}}_\text{error propagation}\\
+\underbrace{\zelfkalkdun{\VAaprealoperator_{\VAtimestep}(f^{(\VAtimestepno)})-\VAapkindifoperator_{\VAtimestep}(f^{(\VAtimestepno)})}}_\text{local error}\,.\label{eq:ap_limnul_errorpartsW1}
\end{multline}
The error at time $t^{(n+1)}$ thus consists of the error at time $t^{(n)}$ that is propagated over the time step, augmented by the local error due to using the approximate KD process instead of the exact kinetic process.

In the remainder of this section, we will use Equation~\eqref{eq:ap_limnul_errorpartsW1} to derive a bound on the Wasserstein error at some final time $\VAtotaltime=\VAnoftimesteps\VAtimestep$ for finite $\VAscaleparam$. We will first treat the last term, expressing the local error due to the difference in distribution when a kinetic process is replaced by a diffusion step, in Section~\ref{subsec:ap_limnul_KD_distrerror}. Then, we consider how the error propagates through different time steps and show the resulting error at time $\VAtotaltime=\VAnoftimesteps\VAtimestep$ in Section~\ref{subsec:ap_limnul_KD_propagation}.

\subsection{Local error when the scaling parameter is finite\label{subsec:ap_limnul_KD_distrerror}}

The KD operator $\VAapkindifoperator_{\VAtimestep}$ is only different from the kinetic operator $\VAaprealoperator_{\VAtimestep}$ when a collision occurs in the time step. Hence, the paths without a collision do not contribute to the term $\zelfkalk{\VAaprealoperator_{\VAtimestep}(f^{(\VAtimestepno)})-\VAapKDoperator_{\VAtimestep}(f^{(\VAtimestepno)})}$ of Equation~\eqref{eq:ap_limnul_errorpartsW1}. We thus only need to consider the paths conditioned on the occurrence of a collision, which we indicate in our operators with a superscript $>0$ as
\begin{align}
\VAaprealoperatorD_{\VAtimestep}(f^{(\VAtimestepno)})&=\VAaprealoperator_{\VAtimestep}(f^{(\VAtimestepno)}|\#\text{collisions in }[t^{(\VAtimestepno)},t^{(\VAtimestepno+1)}]>0)\\
\VAapKDoperatorD_{\VAtimestep}(f^{(\VAtimestepno)})&=\VAapKDoperator_{\VAtimestep}(f^{(\VAtimestepno)}|\#\text{collisions in }[t^{(\VAtimestepno)},t^{(\VAtimestepno+1)}]>0)\,.
\end{align}
The local error term can thus be written as
\begin{multline}
\zelfkalkdun{\VAaprealoperator_{\VAtimestep}(f^{(\VAtimestepno)})-\VAapKDoperator_{\VAtimestep}(f^{(\VAtimestepno)})}\\
=\text{P}\left(\#\text{collisions in }\VAtimestep>0\right)\zelfkalkdun{\VAaprealoperatorD_{\VAtimestep}(f^{(\VAtimestepno)})-\VAapKDoperatorD_{\VAtimestep}(f^{(\VAtimestepno)})}\,.\label{eq:ap_limnul_KD_distrerr_condD_real}
\end{multline}
In the remainder of this section, we will first find an expression for the probability on the right hand side of Equation~\eqref{eq:ap_limnul_KD_distrerr_condD_real} and then bound the second factor, to obtain a bound for the local error. At the end of this section, we will numerically show this is not a sharp bound in general.

\paragraph{Split the local error based on the number of events.} Paths without a collision constitute a fraction $e^{-\frac{\VAratespatial\VAtimestep}{\VAscaleparam^2}}=1-\frac{\VAratespatial\VAtimestep}{\VAscaleparam^2}+\mathcal{O}\left(\VAtimestep^2\right)$, $\VAtimestep\rightarrow0$ of the population with $\mathcal{O}$ the Landau symbol. Hence, in the limit $\VAtimestep\rightarrow0$, we have
\begin{multline}
\zelfkalkdun{\VAaprealoperator_{\VAtimestep}(f^{(\VAtimestepno)})-\VAapKDoperator_{\VAtimestep}(f^{(\VAtimestepno)})}\\
=\frac{\VAratespatial\VAtimestep}{\VAscaleparam^2}\left(1+\mathcal{O}\left(\VAtimestep\right)\right)\zelfkalkdun{\VAaprealoperatorD_{\VAtimestep}(f^{(\VAtimestepno)})-\VAapKDoperatorD_{\VAtimestep}(f^{(\VAtimestepno)})}\,.\label{eq:ap_limnul_KD_distrerr_condD}
\end{multline}
To bound the factor $\zelfkalkdun{\VAaprealoperatorD_{\VAtimestep}(f^{(\VAtimestepno)})-\VAapKDoperatorD_{\VAtimestep}(f^{(\VAtimestepno)})}$ in Equation~\eqref{eq:ap_limnul_KD_distrerr_condD}, we first introduce the additional operator $\VAapkinkinoperatorD_{\VAtimestep}$, indicating a kinetic process with exactly one collision. Via a triangle inequality, we can bound Equation~\eqref{eq:ap_limnul_KD_distrerr_condD} as
\begin{multline}
\zelfkalkdun{\VAaprealoperator_{\VAtimestep}(f(t^{(n)}))-\VAapKDoperator_{\VAtimestep}(f(t^{(n)}))}\\
\leq\frac{\VAratespatial\VAtimestep}{\VAscaleparam^2}\left(1+\mathcal{O}\left(\VAtimestep\right)\right)\zelfkalkdun{\VAaprealoperatorD_{\VAtimestep}(f(t^{(n)}))-\VAapkinkinoperatorD_{\VAtimestep}(f(t^{(n)}))}\\
+\frac{\VAratespatial\VAtimestep}{\VAscaleparam^2}\left(1+\mathcal{O}\left(\VAtimestep\right)\right)\zelfkalkdun{\VAapkinkinoperatorD_{\VAtimestep}(f(t^{(n)}))-\VAapKDoperatorD_{\VAtimestep}(f(t^{(n)}))}\,.\label{eq:ap_limnul_KD_1colsplitup}
\end{multline}

\paragraph{The first term of Equation~\eqref{eq:ap_limnul_KD_1colsplitup}.} $\zelfkalkdun{\VAaprealoperatorD_{\VAtimestep}(f(t^{(n)}))-\VAapkinkinoperatorD_{\VAtimestep}(f(t^{(n)}))}$ captures differences between a kinetic process with at least one collision and a kinetic process with exactly one collision. Only if a second collision takes place, the effect of these operators differs. If the first collision takes place after a time $\tau=\VAtimestep-\theta$, a second collision takes place with a probability $\frac{\VAratespatial\tau}{\VAscaleparam^2}$ and the expected difference in distance in such a case is of order $\theta$ as well, resulting in a contribution proportional to $\VAtimestep\theta^2$. With $\theta\leq\VAtimestep$, we can write
\begin{equation}
\frac{\VAratespatial\VAtimestep}{\VAscaleparam^2}\left(1+\mathcal{O}\left(\VAtimestep\right)\right)\zelfkalkdun{\VAaprealoperatorD_{\VAtimestep}(f(t^{(n)}))-\VAapkinkinoperatorD_{\VAtimestep}(f(t^{(n)}))}\rightarrow\mathcal{O}(\VAtimestep^3),\quad\VAtimestep\rightarrow0\,,
\end{equation}
for the first term of Equation~\eqref{eq:ap_limnul_KD_1colsplitup}. In the remainder of this section, we will derive a bound for the second term of the right hand side of Equation~\eqref{eq:ap_limnul_KD_1colsplitup} which will turn out to be of order $2.5$ in $\VAtimestep$ and will be the dominant term of Equation~\eqref{eq:ap_limnul_KD_1colsplitup}.

\newcommand{\VAstate}{s}
\def\VAtimestepvel{v}
\def\VAtimesteppos{x}
\paragraph{Conditioning to bound the second term of Equation~\eqref{eq:ap_limnul_KD_1colsplitup}.} The second term on the right hand side of Equation~\eqref{eq:ap_limnul_KD_1colsplitup} expresses the difference between a kinetic process with exactly one collision taking place and a KD process where exactly one collision takes place. Until this single collision, both operators behave identically, but after this collision, the kinetic operator $\VAapkinkinoperatorD_{\VAtimestep}$ lets the particle move with a newly sampled velocity whereas the KD operator $\VAapKDoperatorD_{\VAtimestep}$ moves the particle diffusively conditioned on a newly sampled final velocity. To bound the Wasserstein distance between $\VAapkinkinoperatorD_{\VAtimestep}$ and $\VAapKDoperatorD_{\VAtimestep}$, we restrict the joint distribution in Equation~\eqref{eq:ap_intro_W1def} to a subset of all joint distributions $\mathcal{J}(g,h)$ by coupling the paths in two ways. This conditioning can thus lead to overestimating the error.

The first coupling is based on the initial state at the beginning of the timestep $\left(\VAtimesteppos_\VAtimestepno,\VAtimestepvel_\VAtimestepno\right)=\left(x(\VAtimestepno\VAtimestep),v(\VAtimestepno\VAtimestep)\right)$, denoted compactly as $\VAstate_\VAtimestepno$, and the remaining time $\theta_\VAtimestepno=\VAtimestep-\Delta\VAeventtime_{\VAeventnotilltimestep{\VAtimestepno}+1|\VAtimestepno}$, where the time index is as defined in Figure~\ref{fig:ap_eventtime_tijdsstap}. The resulting bound on the error will decrease with decreasing $\VAtimestep$, proving that the KD simulation is consistent with the kinetic model. At the end of this section, we will numerically show the actual error to be significantly lower.

The second coupling strategy is based on the final velocity in the time step, $\nu_\VAtimestepno=\VAtimestepvel_{{\VAtimestepno+1}}$. This coupling is pertinent when we consider the error propagation in Section~\ref{subsec:ap_limnul_KD_propagation}, since a coupling based on significantly different final velocities would have the particles drift apart during a time step and would generally not correspond to the minimizer in Equation~\eqref{eq:ap_intro_W1def}.

We denote the resulting distance for the coupled subset of particles as
\begin{equation}
\zelfkalkdun{\left.\VAapkinkinoperatorD_{\VAtimestep}(f(t^{(n)}))-\VAapKDoperatorD_{\VAtimestep}(f(t^{(n)}))\right|\VAstate_\VAtimestepno,\theta_\VAtimestepno,\nu_\VAtimestepno}\,.
\end{equation}The resulting restriction to a subset of the possible joint probabilities of Equation~\eqref{eq:ap_intro_W1def} bounds the second term of the right hand side of Equation~\eqref{eq:ap_limnul_KD_1colsplitup} as
\begin{align}
\zelfkalkdun{\VAapkinkinoperatorD_{\VAtimestep}(f(t^{(n)}))\!-\!\VAapKDoperatorD_{\VAtimestep}(f(t^{(n)}))}\!\leq\!\mathbb{E}\!\left[\!\zelfkalkdun{\!\left.\VAapkinkinoperatorD_{\VAtimestep}(f(t^{(n)}))\!-\!\VAapKDoperatorD_{\VAtimestep}(f(t^{(n)}))\right|\!\VAstate_\VAtimestepno,\theta_\VAtimestepno,\nu_\VAtimestepno}\!\right]\!\!.\label{eq:ap_KD_limnul_conditioningbound}
\end{align}
This final inequality can also be interpreted as an application of the subadditivity property of the Wasserstein distance \cite[p.~94]{villani2008optimaltransportWasserstein}.

The remainder of this section is devoted to finding an analytical expression for $\mathbb{E}\!\left[\zelfkalkdun{\left.\VAapkinkinoperatorD_{\VAtimestep}(f(t^{(n)}))\!-\!\VAapKDoperatorD_{\VAtimestep}(f(t^{(n)}))\right|\!\VAstate_\VAtimestepno,\theta_\VAtimestepno,\nu_\VAtimestepno}\right]$.

\paragraph{Determining $\zelfkalkdun{\left.\VAapkinkinoperatorD_{\VAtimestep}(f(t^{(n)}))-\VAapKDoperatorD_{\VAtimestep}(f(t^{(n)}))\right|\VAstate_\VAtimestepno,\theta_\VAtimestepno,\nu_\VAtimestepno}$.} Conditioning on $\VAstate_\VAtimestepno$ and $\theta_\VAtimestepno$ also implies conditioning on the position at which the first collision takes place, which is $\VAtimesteppos_\VAtimestepno+\VAtimestepvel_\VAtimestepno(\VAtimestep-\theta_\VAtimestepno)$. Furthermore conditioning on the final velocity $\nu_\VAtimestepno$ means the additional motion $\Delta X$ by the kinetic operator conditioned on only one collision taking place, $\VAapkinkinoperatorD_{\VAtimestep}$, is also fixed and equals $\frac{\nu_\VAtimestepno}{\VAscaleparam}\theta_\VAtimestepno$. The additional motion from the position $\VAtimesteppos_\VAtimestepno+\VAtimestepvel_\VAtimestepno(\VAtimestep-\theta_\VAtimestepno)$ by the KD operator is a normally distributed step with mean and variance according to Equations~\eqref{eq:ap_fixedv0_mean_result} and~\eqref{eq:ap_fixedv0_var}. A comparison between both comes down to a comparison between a dirac delta distribution $\delta_{\frac{\nu_\VAtimestepno}{\VAscaleparam}\theta_\VAtimestepno}$ with probability density function
\begin{equation}
\delta\left(\Delta\VAdiscreteposeend-\frac{\nu_\VAtimestepno}{\VAscaleparam}\theta_\VAtimestepno\right)\label{eq:ap_limnul_KD_dirac}
\end{equation}
and the normal distribution $\mathcal{N}\left(\text{Eq.}\eqref{eq:ap_fixedv0_mean_result},\text{Eq.}\eqref{eq:ap_fixedv0_var}\right)$ with probability density function
\begin{multline}
\frac{1}{\sqrt{\frac{2\pi}{3}\left((\frac{\nu_\VAtimestepno}{\VAscaleparam}-\VAbackgroundspeedmeanflow)^2+\frac{\VAbackgroundspeedtemp}{\VAscaleparam^2}\right)\frac{\VAratespatial\theta_\VAtimestepno^3}{\VAscaleparam^2}\left(1+\mathcal{O}\left(\theta_\VAtimestepno\right)\right)}}e^{-\frac{\left(\Delta\VAdiscreteposeend-\frac{\nu_\VAtimestepno}{\VAscaleparam}\theta_\VAtimestepno+\mathcal{O}(\theta_\VAtimestepno^2)\right)^2}{\frac{2}{3}\left((\frac{\nu_\VAtimestepno}{\VAscaleparam}-\VAbackgroundspeedmeanflow)^2+\frac{\VAbackgroundspeedtemp}{\VAscaleparam^2}\right)\frac{\VAratespatial\theta_\VAtimestepno^3}{\VAscaleparam^2}\left(1+\mathcal{O}\left(\theta_\VAtimestepno\right)\right)}},\ \ \VAtimestep\rightarrow0\,.\label{eq:ap_limnul_KD_normalinlim}
\end{multline}
To establish the Wasserstein distance between the dirac delta of Equation~\eqref{eq:ap_limnul_KD_dirac} and the normal distribution of Equation~\eqref{eq:ap_limnul_KD_normalinlim}, we first bound it via a triangle inequality,
\begin{equation}
\zelfkalkdun{\left.\VAapkinkinoperatorD_{\VAtimestep}(f(t^{(n)}))\!-\!\VAapKDoperatorD_{\VAtimestep}(f(t^{(n)}))\right|\!\VAstate_\VAtimestepno,\theta_\VAtimestepno,\nu_\VAtimestepno}\!=\!\zelfkalkdun{\delta_{\frac{\nu_\VAtimestepno}{\VAscaleparam}\theta_\VAtimestepno}\!\!-\!\mathcal{N}\!\left(\text{Eq.}\eqref{eq:ap_fixedv0_mean_result},\text{Eq.}\eqref{eq:ap_fixedv0_var}\right)}
\end{equation}
\begin{equation}
\hspace{2cm}\leq\zelfkalkdun{\delta_{\frac{\nu_\VAtimestepno}{\VAscaleparam}\theta_\VAtimestepno}\!\!-\!\delta_{\text{Eq.}\eqref{eq:ap_fixedv0_mean_result}}}+\zelfkalkdun{\delta_{\text{Eq.}\eqref{eq:ap_fixedv0_mean_result}}\!\!-\!\mathcal{N}\left(\text{Eq.}\eqref{eq:ap_fixedv0_mean_result},\text{Eq.}\eqref{eq:ap_fixedv0_var}\right)}\,,\label{eq:ap_limnul_KD_triangle2}
\end{equation}
where we use a dirac delta at the mean of Equation~\eqref{eq:ap_limnul_KD_normalinlim} as an additional distribution. The Wasserstein distance between the dirac delta of Equation~\eqref{eq:ap_limnul_KD_dirac} and a dirac delta at the mean of Equation~\eqref{eq:ap_limnul_KD_normalinlim}, the second term on the right hand side of Equation~\eqref{eq:ap_limnul_KD_triangle2}, is of size $\mathcal{O}(\theta_\VAtimestepno^2)$, $\theta_\VAtimestepno\rightarrow0$. The Wasserstein distance of a normal distribution to a dirac delta positioned at its mean simply equals the expected distance of the normal distribution to its mean, or $\sqrt{\frac{2}{\pi}}$ times the standard deviation of the normal distribution. For the distribution of Equation~\eqref{eq:ap_limnul_KD_normalinlim}, the Wasserstein distance thus becomes
\begin{equation}
\zelfkalkdun{\delta_{\text{Eq.}\eqref{eq:ap_fixedv0_mean_result}}\!-\!\mathcal{N}\!\left(\text{Eq.}\eqref{eq:ap_fixedv0_mean_result},\text{Eq.}\eqref{eq:ap_fixedv0_var}\right)}=\!\sqrt{\frac{2}{3\pi}\!\left(\!\left(\frac{\nu_\VAtimestepno}{\VAscaleparam}\!-\!\VAbackgroundspeedmeanflow\right)^2\!\!\!+\!\frac{\VAbackgroundspeedtemp}{\VAscaleparam^2}\!\right)\!\frac{\VAratespatial\theta_\VAtimestepno^3}{\VAscaleparam^2}}\!+\!\mathcal{O}\!\left(\!\sqrt{\theta_\VAtimestepno^5}\right)\!\,.\label{eq:ap_limnul_KD_W1_NvsD}
\end{equation}
Since the $\mathcal{O}\left(\theta^2\right)$ correction due to the first term of Equation~\eqref{eq:ap_limnul_KD_triangle2} does not alter the dominant term in Equation~\eqref{eq:ap_limnul_KD_W1_NvsD}, Equation~\eqref{eq:ap_limnul_KD_W1_NvsD} is an asymptotic bound for $\zelfkalkdun{\left.\VAapkinkinoperatorD_{\VAtimestep}(f(t^{(n)}))\!-\!\VAapKDoperatorD_{\VAtimestep}(f(t^{(n)}))\right|\!\VAstate_\VAtimestepno,\theta_\VAtimestepno,\nu_\VAtimestepno}$.

\paragraph{Determining the right hand side of Equation~\eqref{eq:ap_KD_limnul_conditioningbound}.} Taking the expected value over possible values of $\VAstate_\VAtimestepno$, $\theta_\VAtimestepno$, and $\nu_\VAtimestepno$ of the right hand side term of Equation~\eqref{eq:ap_limnul_KD_W1_NvsD} results in a bound on $\zelfkalkdun{\VAapkinkinoperatorD_{\VAtimestep}(f(t^{(n)}))-\VAapKDoperatorD_{\VAtimestep}(f(t^{(n)}))}$ as stated by Equation~\eqref{eq:ap_KD_limnul_conditioningbound}. We can use the independence of the right hand side of Equation~\eqref{eq:ap_limnul_KD_W1_NvsD} on $\VAstate_\VAtimestepno$, and the independence of $\theta_\VAtimestepno$ and $\nu_\VAtimestepno$ to simplify the calculation of the expected value to
\begin{multline}
\mathbb{E}_{\VAstate_\VAtimestepno,\theta_\VAtimestepno,\nu_\VAtimestepno}\!\!\left[\sqrt{\frac{2}{3\pi}\left(\left(\frac{\nu_\VAtimestepno}{\VAscaleparam}-\VAbackgroundspeedmeanflow\right)^2+\frac{\VAbackgroundspeedtemp}{\VAscaleparam^2}\right)\frac{\VAratespatial\theta_\VAtimestepno^3}{\VAscaleparam^2}}\right]\\
=\frac{2}{3\pi}\mathbb{E}_{\theta_\VAtimestepno}\!\!\left[\sqrt{\left(\frac{\nu_\VAtimestepno}{\VAscaleparam}-\VAbackgroundspeedmeanflow\right)^2+\frac{\VAbackgroundspeedtemp}{\VAscaleparam^2}}\right]\mathbb{E}_{\nu_\VAtimestepno}\!\!\left[\sqrt{\frac{\VAratespatial\theta_\VAtimestepno^3}{\VAscaleparam^2}}\right]\,,\label{eq:ap_limnul_KD_expectovercond}
\end{multline}
where the expectation over $\theta_\VAtimestepno$ is conditioned on at least one collision taking place, resulting in the formula
\begin{align}
\mathbb{E}_{\theta_\VAtimestepno}\!\!\left[\sqrt{\frac{\VAratespatial\theta_\VAtimestepno^3}{\VAscaleparam^2}}\right]&=\frac{\int_0^{\VAtimestep}\sqrt{\frac{\VAratespatial^3(\VAtimestep-\tau)^3}{\VAscaleparam^6}} e^{-\frac{\VAratespatial(\VAtimestep-\tau)}{\VAscaleparam^2}}\text{d}\tau}{1-e^{-\frac{\VAratespatial\VAtimestep}{\VAscaleparam^2}}}\\
&=\frac{\int_0^{\VAtimestep}\sqrt{\frac{\VAratespatial^3(\VAtimestep-\tau)^3}{\VAscaleparam^6}}\left(1+\mathcal{O}(\VAtimestep)\right)\text{d}\theta}{\frac{\VAratespatial\VAtimestep}{\VAscaleparam^2}\left(1+\mathcal{O}(\VAtimestep)\right)},\quad\VAtimestep\rightarrow0,\quad\tau\leq\VAtimestep\\
&=\frac{2}{5}\sqrt{\frac{\VAratespatial\VAtimestep^3}{\VAscaleparam^2}}\left(1+\mathcal{O}(\VAtimestep)\right),\quad\VAtimestep\rightarrow0\,.\label{eq:ap_limnul_KD_expectovertheta_res}
\end{align}
The other expectation in Equation~\eqref{eq:ap_limnul_KD_expectovercond} is over the final velocity $\nu_\VAtimestepno$, which is normally distributed with mean $\VAscaleparam\VAbackgroundspeedmeanflow$ and variance $\VAbackgroundspeedtemp$. Using the definition of the expectation and a simple change of variable, $\upsilon=\frac{\nu_\VAtimestepno-\VAscaleparam\VAbackgroundspeedmeanflow}{\sqrt{\VAbackgroundspeedtemp}}$, yields
\begin{align}
\mathbb{E}_{\nu_\VAtimestepno}\!\!\left[\sqrt{\left(\frac{\nu_\VAtimestepno}{\VAscaleparam}-\VAbackgroundspeedmeanflow\right)^2+\frac{\VAbackgroundspeedtemp}{\VAscaleparam^2}}\right]&=\int\sqrt{\left(\frac{\nu_\VAtimestepno}{\VAscaleparam}-\VAbackgroundspeedmeanflow\right)^2+\frac{\VAbackgroundspeedtemp}{\VAscaleparam^2}}\frac{1}{\sqrt{2\pi\VAbackgroundspeedtemp}}e^{-\frac{(\nu_\VAtimestepno-\VAscaleparam\VAbackgroundspeedmeanflow)^2}{2\VAbackgroundspeedtemp}}\text{d}\nu_\VAtimestepno\\
&=\sqrt{\frac{\VAbackgroundspeedtemp}{\VAscaleparam^2}}\underbrace{\int\sqrt{\upsilon^2+1}\frac{1}{\sqrt{2\pi}}e^{-\frac{\upsilon^2}{2}}\text{d}\upsilon}_{\approx 1.3545
\text{ (numerically)}}\,.\label{eq:ap_limnul_KD_expectovernu_simplifiedintegral}
\end{align}
The remaining integral in Equation~\eqref{eq:ap_limnul_KD_expectovernu_simplifiedintegral} can be numerically computed to equal $1.3545$. Substituting the results from Equations~\eqref{eq:ap_limnul_KD_expectovertheta_res} and~\eqref{eq:ap_limnul_KD_expectovernu_simplifiedintegral} in Equation~\eqref{eq:ap_limnul_KD_expectovercond}, gives, according to Equation~\eqref{eq:ap_KD_limnul_conditioningbound}, the asymptotic bound 
\begin{equation}
\zelfkalkdun{\VAapkinkinoperatorD_{\VAtimestep}(f(t^{(n)}))\!-\!\VAapKDoperatorD_{\VAtimestep}(f(t^{(n)}))}\leq0.24959\sqrt{\frac{\VAbackgroundspeedtemp\VAratespatial\VAtimestep^3}{\VAscaleparam^4}}+\mathcal{O}\left(\sqrt{\VAtimestep^5}\right),\quad\VAtimestep\rightarrow0\,.\label{eq:ap_limnul_KD_distrerr_res}
\end{equation}

\paragraph{Resulting bound.} As stated before, the bound on $\zelfkalkdun{\VAapkinkinoperatorD_{\VAtimestep}(f(t^{(n)}))\!-\!\VAapKDoperatorD_{\VAtimestep}(f(t^{(n)}))}$ dominates the right hand side of Equation~\eqref{eq:ap_limnul_KD_1colsplitup}, since the first term is at least third order in $\VAtimestep$. We thus find as the bound for the local error term of Equation~\eqref{eq:ap_limnul_errorpartsW1}
\begin{equation}
\zelfkalkdun{\VAaprealoperator_{\VAtimestep}(f(t^{(n)}))-\VAapKDoperator_{\VAtimestep}(f(t^{(n)}))}\leq0.24959\sqrt{\frac{\VAbackgroundspeedtemp\VAratespatial^3\VAtimestep^5}{\VAscaleparam^8}}+\mathcal{O}\left(\sqrt{\VAtimestep^7}\right),\quad\VAtimestep\rightarrow0\,.
\label{eq:ap_limnul_KD_distrerr_res}
\end{equation}

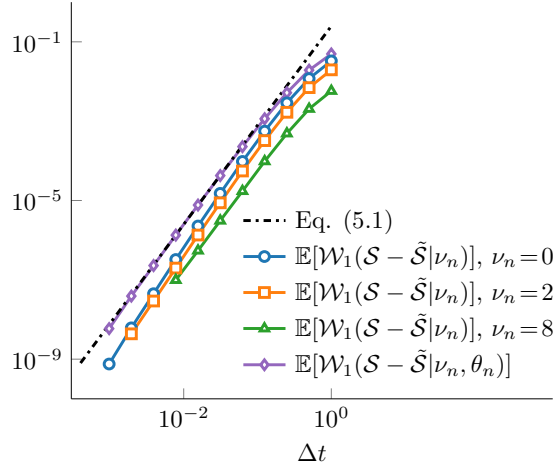
\begin{figure}[H]
\centering
\def\figfolderhomo{figuren/}
\def\zelfhomofigboundA{\ref{eq:ap_limnul_KD_distrerr_res}}
\def\zelfhomofigboundB{\ref{eq:ap_liminf_KD_distrerror}}
\begin{minipage}{.65\textwidth}
  \centering
  \resizebox{\textwidth}{!}{%
	\definecolor{colorA}{RGB}{31,119,180}
\definecolor{colorB}{RGB}{255,127,14}
\definecolor{colorC}{RGB}{44,160,44}
\definecolor{colorD}{RGB}{148,103,189}
\begin{tikzpicture}
\begin{axis}[
  width=240pt,
  height=207pt,
  xmin=3e-4, xmax=1e3,
  ymin=1e-10, ymax=1e0,
  axis y line*=left,
  ytick = {.000000001, .00001, .1},
  xtick = {.0001,.01,1},
  axis x line*=bottom,
  xmode=log,
  ymode=log,
  xlabel=$\VAtimestep$,
  ylabel=$\vphantom{\mathbb{E}[\mathcal{W}_1(\VAaprealoperator-\VAapKDoperator|\VAfinalV,\tau)]/\sqrt{\VAbackgroundspeedtemp}}$,
  xlabel near ticks,
  legend cell align={left},
  legend style={draw=none, fill=none,at={(axis cs:10,0.000002)}}]
\addplot[very thick, color=colorD, mark options={fill=white}, mark=diamond*] table [x index=0, y expr=\thisrowno{1}*(1-exp(-\thisrowno{0}))]{\figfolderhomo KD_lowcol_aggreg2.out};
\label{fig_homo_KD_conform}
\addplot[very thick, color=colorA, mark options={fill=white}, mark=*] table [x index=0, y expr=\thisrowno{1}*(1-exp(-\thisrowno{0}))]{\figfolderhomo KD_lowcol_aggreg_condnu.out};
\label{fig_homo_KD_nocondtheta0}
\addplot[very thick, color=colorB, mark options={fill=white}, mark=square*] table [x index=0, y expr=\thisrowno{1}*(1-exp(-\thisrowno{0})), skip coords between index={10}{11}]{\figfolderhomo KD_lowcol_aggreg_condnu_2.out};
\label{fig_homo_KD_nocondtheta2}
\addplot[very thick, color=colorC, mark options={fill=white}, mark=triangle*] table [x index=0, y expr=\thisrowno{1}*(1-exp(-\thisrowno{0})), skip coords between index={8}{11}]{\figfolderhomo KD_lowcol_aggreg_condnu_8.out};
\label{fig_homo_KD_nocondtheta8}
\addplot[very thick, dashdotted, black, domain=.0004:1,samples=2] {.2496*x^(5/2)};
\label{fig_homo_KD_goodorder_bound}
\node [draw=none,fill=none,anchor=south] at (axis cs: 8,0.0000000002) {\shortstack[l]{
\ref{fig_homo_KD_goodorder_bound} Eq. (\zelfhomofigboundA)\\
\ref{fig_homo_KD_nocondtheta0} $\mathbb{E}[\mathcal{W}_1(\VAaprealoperator-\VAapKDoperator|\VAfinalV_\VAtimestepno)]$, $\VAfinalV_\VAtimestepno\!=\!0$\\
\ref{fig_homo_KD_nocondtheta2} $\mathbb{E}[\mathcal{W}_1(\VAaprealoperator-\VAapKDoperator|\VAfinalV_\VAtimestepno)]$, $\VAfinalV_\VAtimestepno\!=\!2$\\
\ref{fig_homo_KD_nocondtheta8} $\mathbb{E}[\mathcal{W}_1(\VAaprealoperator-\VAapKDoperator|\VAfinalV_\VAtimestepno)]$, $\VAfinalV_\VAtimestepno\!=\!8$\\
\ref{fig_homo_KD_conform} $\mathbb{E}[\mathcal{W}_1(\VAaprealoperator-\VAapKDoperator|\VAfinalV_\VAtimestepno,\theta_\VAtimestepno)]$}};
\end{axis}
\end{tikzpicture}%
}
\end{minipage}
  \captionof{figure}{The Wasserstein distance between the distribution of the kinetic process and of the KD process after conditioning on the final velocity and the time of the diffusion process as a function of $\VAtimestep$, conducted with 200,000 particle paths with at least one collision.}
  \label{fig:ap_num_nt1_KD_cond}
\end{figure}

\paragraph{Numerical validation.} We have validated the bound of Equation~\eqref{eq:ap_limnul_KD_distrerr_res} on the new error during a time step numerically. To do so, we assumed all particles have an identical initial state $\VAstate_\VAtimestepno=(\VAtimesteppos_\VAtimestepno,\VAtimestepvel_\VAtimestepno)$ and we compare the outcome of the kinetic and KD simulation conditioned on the final velocity $\VAfinalV_\VAtimestepno$, with $\VAratespatial=\VAscaleparam=\VAbackgroundspeedtemp=1$ and $\VAbackgroundspeedmeanflow=1$. Then, for different $\VAtimestepvel_\VAtimestepno$, we see in Figure~\ref{fig:ap_num_nt1_KD_cond} that the analytical bound of Equation~\eqref{eq:ap_limnul_KD_distrerr_res} overestimates the actual error. This difference originates from the conditioning on the remaining time $\theta_\VAtimestepno$, which was also included in the analysis. If we do include the conditioning on $\theta_\VAtimestepno$, we achieve a tight fit with the overestimating bound, as also shown in Figure~\ref{fig:ap_num_nt1_KD_cond}. We thus conclude the error bound of Equation~\eqref{eq:ap_limnul_KD_distrerr_res} to be valid, but not tight. In the next section, we show the error at a final time $\VAtotaltime=\VAnoftimesteps\VAtimestep$ based on this error bound, and conclude convergence to the kinetic process as $\VAtimestep\rightarrow0$.

\subsection{Error propagation and total error when the scaling parameter is finite\label{subsec:ap_limnul_KD_propagation}}

In this Section, we first find a bound on the error akin to Equation~\eqref{eq:ap_limnul_errorpartsW1}, but with conditioning on the velocity at the end of the time steps. This conditioning results in an overestimation of the error, but it enables the computation of the error propagation from time $t^{(n)}$ to time $t^{(n+1)}$. Then, we combine the result with the outcome of the local error analysis in Section~\ref{subsec:ap_limnul_KD_distrerror} to find the total error at the simulation end time $\VAtotaltime=\VAnoftimesteps\VAtimestep$.

Since Equation~\eqref{eq:ap_limnul_errorpartsW1} also holds when conditioning on the velocities $\{\VAtimestepvel_\VAtimestepno\}_{\VAtimestepno=0}^{\VAnoftimesteps}=\{v(\VAtimestepno\VAtimestep)\}_{\VAtimestepno=0}^{\VAnoftimesteps}$ at the discrete times $\{\VAtimestepno\VAtimestep\}_{\VAtimestepno=0}^{\VAnoftimesteps}$, we can write
\begin{multline}
\mathbb{E}\left[\zelfkalkdun{\left.f(t^{(\VAtimestepno+1)})-f^{(\VAtimestepno+1)}\right|\{\VAtimestepvel_{{\VAtimestepno}}\}_{\VAtimestepno=0}^{\VAnoftimesteps}}\right]\\
\leq\underbrace{\mathbb{E}\left[\zelfkalkdun{\left.f(t^{(\VAtimestepno)})+\VAaprealoperator_{\VAtimestep}(f(t^{(\VAtimestepno)}))-\left(f^{(\VAtimestepno)}+\VAaprealoperator_{\VAtimestep}(f^{(\VAtimestepno)})\right)\right|\{\VAtimestepvel_{{\VAtimestepno}}\}_{\VAtimestepno=0}^{\VAnoftimesteps}}\right]}_\text{error propagation}\\
+\underbrace{\mathbb{E}\left[\zelfkalkdun{\left.\VAaprealoperator_{\VAtimestep}(f^{(\VAtimestepno)})-\VAapkindifoperator_{\VAtimestep}(f^{(\VAtimestepno)})\right|\{\VAtimestepvel_{{\VAtimestepno}}\}_{\VAtimestepno=0}^\VAnoftimesteps}\right]}_\text{local error}\,,\label{eq:ap_limnul_errorpartsW1_condVs}
\end{multline}
and we note that the subadditivity property
\begin{equation}
\zelfkalkdun{f(t^{(\VAtimestepno)})-f^{(\VAtimestepno)}}\leq\mathbb{E}\left[\zelfkalkdun{\left.f(t^{(\VAtimestepno)})-f^{(\VAtimestepno)}\right|\{\VAtimestepvel_{{\VAtimestepno}}\}_{\VAtimestepno=0}^{\VAnoftimesteps}}\right]
\end{equation}
holds for all $\VAtimestepno$.

The error propagation term from Equation~\eqref{eq:ap_limnul_errorpartsW1_condVs} can easily be bounded by the above conditioning on the velocities, since the kinetic operator will have an identical effect on the population that has the same initial and final velocity, and thus
\begin{multline}
\mathbb{E}\left[\zelfkalkdun{\left.f(t^{(\VAtimestepno)})+\VAaprealoperator_{\VAtimestep}(f(t^{(\VAtimestepno)})-\left(f^{(\VAtimestepno)}+\VAaprealoperator_{\VAtimestep}(f^{(\VAtimestepno)})\right)\right|\{\VAtimestepvel_{{\VAtimestepno}}\}_{\VAtimestepno=0}^{\VAnoftimesteps}}\right]\\
\leq\mathbb{E}\left[\zelfkalkdun{\left.f(t^{(\VAtimestepno)})-f^{(\VAtimestepno)}\right|\{\VAtimestepvel_{{\VAtimestepno}}\}_{\VAtimestepno=0}^{\VAnoftimesteps}}\right].
\end{multline}
Since we also derived a bound for the local error after conditioning on the velocity at time $t^{(\VAtimestepno)}=\VAtimestepno\VAtimestep$ in Section~\ref{subsec:ap_limnul_KD_distrerror} by conditioning on the final velocity and on the assumption of an identical initial state, we can use the bound from Equation~\eqref{eq:ap_limnul_KD_distrerr_res} as a bound for the local error term in Equation~\eqref{eq:ap_limnul_errorpartsW1_condVs}.

Combining the above elements means we can bound the right hand side of Equation~\eqref{eq:ap_limnul_errorpartsW1_condVs} as
\begin{multline}
\mathbb{E}\left[\zelfkalkdun{\left.f(t^{(\VAtimestepno+1)})-f^{(\VAtimestepno+1)}\right|\{v_\VAtimestepno\}_{\VAtimestepno=0}^\VAnoftimesteps}\right]
\leq\mathbb{E}\left[\zelfkalkdun{\left.f(t^{(\VAtimestepno)})-f^{(\VAtimestepno)}\right|\{v_\VAtimestepno\}_{\VAtimestepno=0}^\VAnoftimesteps}\right]\\
+0.24959\sqrt{\frac{\VAbackgroundspeedtemp\VAratespatial^3\VAtimestep^5}{\VAscaleparam^{8}}}+\mathcal{O}\left(\VAtimestep^3\right)\,,\qquad\VAtimestep\rightarrow0\,.
\end{multline}
With an initial error of zero, this recursive equation has as a solution at time $\VAtotaltime=\VAnoftimesteps\VAtimestep$ which bounds the error at that final time:
\begin{equation}
\zelfkalkdun{f(t^{(\VAtimestepno+1)})-f^{(\VAtimestepno+1)}}\leq0.24959\VAtotaltime\sqrt{\frac{\VAbackgroundspeedtemp\VAratespatial^3\VAtimestep^3}{\VAscaleparam^{8}}}+\mathcal{O}\left(\VAtimestep^2\right),\qquad\VAtimestep\rightarrow0\,.
\end{equation}
This proves the new KD simulation scheme becomes identical to a standard, unbiased, particle tracing simulation of the Boltzmann-BGK equation when $\VAtimestep\rightarrow0$.

\section{Asymptotic preserving property and convergence analysis in high-collisional regimes\label{sec:ap_liminf}}

In this section, we consider the error in the diffusive regime, i.e., when $\VAscaleparam$ becomes very small. In contrast to the treatment of in Section~\ref{sec:ap_limnul}, the correlation between subsequent time steps is not a dominant feature in this regime. The first kinetic flight in a time step, which is the origin of this correlation, only makes up an expected fraction $\mathcal{O}\left(\VAscaleparam^2\right)$, $\VAscaleparam\rightarrow0$ of the time step. Furthermore, the execution of this initial flight is identical in the kinetic simulation and the KD simulation. As a consequence, we can ignore the initial flight. We thus consider the kinetic simulation beginning with a resampled velocity, which is represented by the newly introduced operator $\VAapKoperatoronlyD_{\VAtimestep}$, and we consider the KD scheme only consisting of a diffusion part, which is denoted by the newly introduced operator $\VAapDoperator_{\VAtimestep}$. As we just argued, these operators become identical to the actual kinetic and KD operators in the diffusive limit, i.e.,
\begin{align}
\VAapKoperatoronlyD_{\VAtimestep}\rightarrow\VAaprealoperator_{\VAtimestep}\text{  and  }\VAapDoperator_{\VAtimestep}\rightarrow\VAapKDoperator_{\VAtimestep}\,,\quad\VAscaleparam\rightarrow0\,.
\end{align}
We use the above defined operators to again write the Wasserstein error at time $t^{(\VAtimestepno+1)}$ in terms of the error propagated from time $t^{(\VAtimestepno)}$ and a new local error
\begin{multline}
\zelfkalkdun{f(t^{(\VAtimestepno+1)})\!-\!f^{(\VAtimestepno+1)}}\leq\underbrace{\zelfkalkdun{f(t^{(\VAtimestepno)})\!+\!\VAapKoperatoronlyD_{\VAtimestep}(f(t^{(\VAtimestepno)}))\!-\!\left(f^{(\VAtimestepno)}
\!+\!\VAapKoperatoronlyD_{\VAtimestep}(f^{(\VAtimestepno)})\right)}}_\text{error propagation}\\
+\underbrace{\zelfkalkdun{\VAapKoperatoronlyD_{\VAtimestep}(f^{(\VAtimestepno)})-\VAapDoperator_{\VAtimestep}(f^{(\VAtimestepno)})}}_\text{local error}\,,\ \ \VAscaleparam\rightarrow0\,,\label{eq:ap_liminf_errorpartsW1}
\end{multline}
where the usage of the operators without the inital kinetic flight forms an additional overestimation of the error. The derivation will follow the same structure as in the low-collisional case. We will first find the local error during a time step in Section~\ref{subsec:ap_liminf_KD_distrerror}. Then, we discuss the error propagation and find the total error bound via a geometric series in Section~\ref{subsec:ap_liminf_KD_propagation}.

\subsection{Local error in the diffusive limit\label{subsec:ap_liminf_KD_distrerror}}

The final term in Equation~\eqref{eq:ap_liminf_errorpartsW1} expresses the difference between a diffusion step and a kinetic step starting from the same initial condition $f^{(\VAtimestepno)}$. We bound this local error with the subadditivity property of the Wasserstein distance via conditioning on the initial position $x(\VAtimestepno\VAtimestep)=x_\VAtimestepno$ at time $t^{(\VAtimestepno)}$ and on the final velocity $\nu_\VAtimestepno$ in a similar fashion as in Equation~\eqref{eq:ap_KD_limnul_conditioningbound},
\begin{equation}
\zelfkalkdun{\VAapKoperatoronlyD_{\VAtimestep}(f^{(\VAtimestepno)})-\VAapDoperator_{\VAtimestep}(f^{(\VAtimestepno)})}\leq\mathbb{E}\left[\zelfkalkdun{\left.\VAapKoperatoronlyD_{\VAtimestep}(f^{(\VAtimestepno)})-\VAapDoperator_{\VAtimestep}(f^{(\VAtimestepno)})\right|x_\VAtimestepno,\nu_\VAtimestepno}\right]\,.\label{eq:ap_liminf_KD_conditioningbound}
\end{equation}
The above conditioning is equivalent to having a dirac delta at $x_\VAtimestepno$, $\delta_{x_\VAtimestepno}$, as the initial distribution $f^{(\VAtimestepno)}$. Note that the velocity information at $f^{(\VAtimestepno)}$ is not used by either of the operators $\VAapKoperatoronlyD_{\VAtimestep}$ and $\VAapDoperator_{\VAtimestep}$, because they consider velocity resampling and only a Brownian increment, respectively. We can thus write
\begin{equation}
\zelfkalkdun{\left.\VAapKoperatoronlyD_{\VAtimestep}(f^{(\VAtimestepno)})-\VAapDoperator_{\VAtimestep}(f^{(\VAtimestepno)})\right|x_\VAtimestepno,\nu_\VAtimestepno}=\zelfkalkdun{\left.\VAapKoperatoronlyD_{\VAtimestep}(\delta_{x_\VAtimestepno})-\VAapDoperator_{\VAtimestep}(\delta_{x_\VAtimestepno})\right|x_\VAtimestepno,\nu_\VAtimestepno}\,.
\end{equation}

To obtain an estimate of $\zelfkalk{\left.\VAapKoperatoronlyD_{\VAtimestep}\left(\delta_{\VAtimesteppos_\VAtimestepno}\right)-\VAapDoperator_{\VAtimestep}\left(\delta_{\VAtimesteppos_\VAtimestepno}\right)\right|x_\VAtimestepno,\nu_\VAtimestepno}$ in the diffusive limit, we perform an Edgeworth expansion~\cite{hall2013edgeworth} of the distribution $\delta_{\VAtimesteppos_\VAtimestepno}+\VAapKoperatoronlyD_{\VAtimestep}\left(\delta_{\VAtimesteppos_\VAtimestepno}\right)$ that is the result of kinetic evolution of the initial condition $\delta_{\VAtimesteppos_\VAtimestepno}$ over a time step $\VAtimestep$. The Edgeworth expansion approximates a probability distribution in terms of its cumulants. The cumulants of a random variable $X$, here the position at time $(\VAtimestepno+1)\VAtimestep$, are defined via the cumulant-generating function $\mathcal{K}'(t)$,
\begin{equation}
\mathcal{K}'(\VAtimestep)=\log\left(\mathbb{E}\left[e^{\VAtimestep X}\right]\right)=\sum_{i=1}^\infty\kappa_i\frac{\VAtimestep^i}{i!}=m\VAtimestep+\frac{s^2\VAtimestep^2}{2}+\sum_{i=3}^\infty\kappa_i\frac{\VAtimestep^i}{i!}\,,
\end{equation}
in which the expectation is taken with respect to the probability distribution $\delta_{\VAtimesteppos_\VAtimestepno}+\VAapKoperatoronlyD_{\VAtimestep}\left(\delta_{\VAtimesteppos_\VAtimestepno}\right)$, $\kappa_i$ denotes the $i$-th cumulant, and $m=\kappa_1$ and $s^2=\kappa_2$ are the mean and variance of $\delta_{\VAtimesteppos_\VAtimestepno}+\VAapKoperatoronlyD_{\VAtimestep}\left(\delta_{\VAtimesteppos_\VAtimestepno}\right)$. By construction, the position distribution of the diffusion process, denoted as $\delta_{\VAtimesteppos_\VAtimestepno}+\VAapDoperator_{\VAtimestep}\left(\delta_{\VAtimesteppos_\VAtimestepno}\right)$ is a normal distribution with the same mean and variance as that of the kinetic process, meaning the cumulant generating function is
\begin{equation}
\hat{\mathcal{K}}(\VAtimestep)=m\VAtimestep+\frac{s^2\VAtimestep^2}{2}\,.
\end{equation}
We can therefore write the Edgeworth expansion of the distribution
\begin{equation}
\delta_{x_\VAtimestepno}(x)\!+\!\VAapKoperatoronlyD_{\VAtimestep}(\delta_{\VAtimesteppos_\VAtimestepno})(x)\!=\!\left(\delta_{x_0}(x)\!+\!\VAapDoperator_{\VAtimestep}(\delta_{\VAtimesteppos_\VAtimestepno})(x)\right)\!\!\left(\!1+\sum_{i=3}^\infty\frac{\kappa_i}{i!\VAewstdv^i}\VAhermite{i}\!\left(\!\frac{x\!-\!\VAewmean}{\VAewstdv}\!\right)\!\right)\,,\label{eq:ap_KD_liminf_homo_EW}
\end{equation}
in which $\VAhermite{i}$ denotes the Hermite polynomial of degree $i$. Bounding the Wasserstein distance $\zelfkalk{\left.\VAapKoperatoronlyD_{\VAtimestep}\left(\delta_{\VAtimesteppos_\VAtimestepno}\right)-\VAapDoperator_{\VAtimestep}\left(\delta_{\VAtimesteppos_\VAtimestepno}\right)\right|x_\VAtimestepno,\nu_\VAtimestepno}$ therefore reduces to bounding the cumulants $\kappa_i$.

\paragraph{Bounding the odd cumulants.} To obtain an estimate on the odd cumulants, it is sufficient to realize that the dominant odd cumulant can be expressed in terms of odd standardized moments~\cite{hall2013edgeworth}, on which we can reason more intuitively. The standardized moments $m_i$ are defined as
\begin{equation}
m_i=\mathbb{E}\left[\left(\frac{\Delta X-m}{s}\right)^i\right]\,.
\end{equation}
When $\VAscaleparam\rightarrow0$, the dominant odd term of Equation~\eqref{eq:ap_KD_liminf_homo_EW} will be formed by the third moment $m_3$ of the distribution $\delta_{x_\VAtimestepno}(x)+\VAapKoperatoronlyD_{\VAtimestep}\left(\delta_{\VAtimesteppos_\VAtimestepno}\right)\left(x\right)$ of the position $X$ due to the kinetic process conditioned on the initial position $\VAtimesteppos_\VAtimestepno$ and on the final velocity $\VAdiscreteveleend_\VAnofeventsintimestepnon$. The increment $\Delta X=X-x_\VAtimestepno$ arises as the sum of flight path contributions, as expressed by Equation~\eqref{eq:ap_deltax_assum}. Following on the computation of the first and second moment of the $\Delta X$ in Section~\ref{subsec:ap_fixedv0}, the third moment can also be found. Here, we will however only point out that the third moment arises due to asymmetry of the distribution of $\Delta X$ and the only asymmetry in Equation~\eqref{eq:ap_deltax_assum} arises due to the final flight $\frac{\VAdiscreteveleend_\VAnofeventsintimestepnon}{\VAscaleparam}\VAeventtimedif_\VAnofeventsintimestepnon$. With the expected duration of the final flight expressed by Equation~\eqref{eq:ap_fixedv0_dt0mean}, the third central moment of this final flight conditioned on the final velocity is found to be
\begin{equation}
\mathbb{E}\left[\left(\frac{\VAdiscreteveleend_\VAnofeventsintimestepnon}{\VAscaleparam}\VAeventtimedif_\VAnofeventsintimestepnon-\frac{\VAdiscreteveleend_\VAnofeventsintimestepnon}{\VAscaleparam}\frac{\VAscaleparam^2}{\VAratespatial}\left(1-e^{-\frac{\VAratespatial\VAtimestep}{\VAscaleparam^2}}\right)\right)^3\right]=\mathcal{O}\left(\VAscaleparam^3\right),\ \ \VAscaleparam\rightarrow0\,.
\end{equation}
This gives the size of $m_3$, since $s\rightarrow\mathcal{O}(1)$, $\VAscaleparam\rightarrow0$. In the next paragraph, we will see the dominant even error term of Equation~\eqref{eq:ap_KD_liminf_homo_EW} is asymptotically larger, since it is proportional to $\VAscaleparam^2$, $\VAscaleparam\rightarrow0$.

\paragraph{Bounding the even cumulants.} For the even cumulants, we take a different viewpoint by considering the time step $\VAtimestep$ as consisting of $\frac{\VAratespatial\VAtimestep}{\VAscaleparam^2}$ subtimesteps of equal duration $\frac{\VAscaleparam^2}{\VAratespatial}$. We thus split the positional increment $\Delta X$ over its contributions during these subtimesteps as
\begin{equation}
\Delta X=\sum_{j=1}^{\VAratespatial\VAtimestep/\VAscaleparam^2}\Delta X_j\,.
\end{equation}
Except for an overall negligible effect due to the final flight path, the different $\Delta X_j$ are identically distributed but not independent, since subtimesteps are correlated if no collision takes place during intermediate subtimesteps. This correlation will only be significant for a few neighbouring subtimesteps since, for the $j$-th and $j'$-th subtimesteps, the probability of no collision taking place in the $|j'-j|-1$ intermediate subtimesteps, is
\begin{equation}
e^{-\frac{\VAratespatial}{\VAscaleparam^2}(|j'-j|-1)\frac{\VAscaleparam^2}{\VAratespatial}}=e^{1-|j'-j|}\,,
\end{equation}
which is independent of $\VAscaleparam$ and decreases exponentially as a function of the number of intermediate subtimesteps. If the $\Delta X_j$ would be independent, the cumulant for $\Delta X$ would equal $\frac{\VAratespatial\VAtimestep}{\VAscaleparam^2}$ times the cumulant of a single $\Delta X_j$. For all $\Delta X_j$, the correlation with other subtimesteps is identical, except for those close to the beginning and end of the full time step, since fewer correlated subtimesteps are present for those. When $\VAscaleparam\rightarrow0$, and the number of subtimesteps in which the time step is split thus increases, the portion of subtimesteps at the edges decreases as $\VAscaleparam^2$ and we can effectively consider the time step as a combination of a number of independent, identically distributed subtimesteps. The correlation between different $\Delta X_j$ means each increment has a larger effective contribution than if it is considered in isolation, but there is a lower effective number of subtimesteps. We can capture these effects in a multiplicative constant $k'_i$ to obtain
\begin{equation}
\kappa_i\rightarrow k'_i\frac{\VAratespatial\VAtimestep}{\VAscaleparam^2}\tilde{\kappa}_i,\ \ \VAscaleparam\rightarrow0,\ \ \text{ if }i\text{ is even}\,,
\end{equation}
with $\tilde{\kappa}_i$ the cumulant of the kinetic process over a time $\frac{\VAscaleparam^2}{\VAratespatial}$. The cumulant $\tilde{\kappa}_i$ of a single part equals a constant times $\left(\frac{\sqrt{\VAscaleparam^2\VAbackgroundspeedtemp}}{\VAratespatial}\right)^i$, since the deviation of the mean is a constant times $\frac{\sqrt{\VAbackgroundspeedtemp}}{\VAscaleparam}$ and the size of the step is $\frac{\VAscaleparam^2}{\VAratespatial}$ when $\frac{\VAratespatial\VAtimestep}{\VAscaleparam^2}$ is high. In conclusion,
\begin{equation}
\kappa_i\rightarrow k_i\frac{\sqrt{\VAbackgroundspeedtemp^i}\VAscaleparam^{i-2}}{\VAratespatial^{i-1}\VAtimestep^{i-1}},\ \ \VAscaleparam\rightarrow0,\ \ \text{ if }i\text{ is even}\,,\label{eq:ap_KD_liminf_homo_cumeven}
\end{equation}
with $k_i$ a constant.

\paragraph{Determining the dominant term.} The dominant error term is the term with $i=4$ in Equation~\eqref{eq:ap_KD_liminf_homo_EW}, since it results in a second power of $\VAscaleparam$ as can be seen from Equation~\eqref{eq:ap_KD_liminf_homo_cumeven}. For this dominant term, we use the limiting value for $\VAewstdv$, which equals $\sqrt{2\frac{\VAbackgroundspeedtemp}{\VAratespatial}\VAtimestep}$ as given at the end of section~\ref{subsec:ap_fixedv0}. The dominant error contribution due to a difference in distribution thus equals
{\def\zelflokaalhoogste{\VAapKoperatoronlyD_{\VAtimestep}(f^{(\VAtimestepno)})\!-\!\VAapDoperator_{\VAtimestep}(f^{(\VAtimestepno)})}\begin{align}
\zelfkalkdun{\left.\zelflokaalhoogste\right.\right.&\left.\left.\!\!\!\!\vphantom{\zelflokaalhoogste}\right|\VAtimesteppos_\VAtimestepno,\nu_\VAtimestepno}=\zelfkalk{\frac{1}{\sqrt{2\pi}\VAewstdv}\!e^{-\frac{(x-\mu)^2}{2\VAewstdv^2}}\!\!\!\!\frac{k_4\VAscaleparam^2}{4!\sqrt{8}\VAratespatial\VAtimestep}\VAhermite{4}\!\left(\!\frac{x-\mu}{\VAewstdv}\!\right)\!}\\
&=\frac{k_4\VAscaleparam^2\VAewstdv}{4!\sqrt{8}\VAratespatial\VAtimestep}\zelfkalk{\frac{1}{\sqrt{2\pi}}e^{-\frac{x^2}{2}}(x^4-6x^2+3)}\\
&=\frac{k_4\sqrt{\VAbackgroundspeedtemp}\VAscaleparam^{3}}{4!4\sqrt{\VAratespatial^{3}\VAtimestep}}\zelfkalk{\frac{1}{\sqrt{2\pi}}e^{-\frac{x^2}{2}}(x^4-6x^2+3)}\,,\label{eq:ap_liminf_KD_distrerror_domterm}
\end{align}}
as $\frac{\VAratespatial\VAtimestep}{\VAscaleparam^2}\rightarrow\infty$. The values for $k_4$ and $\zelfkalk{\frac{1}{\sqrt{2\pi}}e^{-x^2/2}(x^4-6x^2+3)}$ have been computed numerically, resulting in values $18.3$ and $1.51$. 
Since neither $\VAtimesteppos_\VAtimestepno$ nor $\nu_\VAtimestepno$ feature in the right hand side of Equation~\eqref{eq:ap_liminf_KD_distrerror_domterm}, taking the expectation over all values of $\VAtimesteppos_\VAtimestepno$ and $\nu_\VAtimestepno$ does not alter the result. In accordance to Equation~\eqref{eq:ap_liminf_KD_conditioningbound}, we thus find the bound
\begin{equation}
\zelfkalk{\VAapKoperatoronlyD_{\VAtimestep}(f^{(\VAtimestepno)})-\VAapDoperator_{\VAtimestep}(f^{(\VAtimestepno)})}\leq 0.58\sqrt{\VAbackgroundspeedtemp}\frac{\VAscaleparam^{3}}{\sqrt{\VAratespatial^{3}\VAtimestep}}\label{eq:ap_liminf_KD_distrerror}\,.
\end{equation}

\begin{figure}[H]
\centering
\def\figfolderhomo{figuren/}
\def\zelfhomofigboundA{\ref{eq:ap_limnul_KD_distrerr_res}}
\def\zelfhomofigboundB{\ref{eq:ap_liminf_KD_distrerror}}
\begin{minipage}{.65\textwidth}
  \centering
  \resizebox{\textwidth}{!}{%
	\definecolor{colorA}{RGB}{31,119,180}
\definecolor{colorB}{RGB}{255,127,14}
\definecolor{colorC}{RGB}{44,160,44}
\definecolor{colorD}{RGB}{148,103,189}
\begin{tikzpicture}
\begin{axis}[
  width=240pt,
  height=207pt,
  ymin=.75e-4, ymax=3e-2,
  xmin=.25e-2, xmax=2.5e2,
  axis y line*=left,
  ytick = {.0001, .001, .01},
  xtick = {1,10,100,1000},
  xminorticks=false,
  yminorticks=false,
  axis x line*=bottom,
  xmode=log,
  ymode=log,
  xlabel=$\VAratespatial/\VAscaleparam^2$,
  ylabel=$\mathcal{W}_1$ error,
  xlabel near ticks,
  legend cell align={left},
  legend style={draw=none, fill=none, at={(axis cs:7,.000002)}, anchor=south}]
\addplot[very thick, color=colorA, mark options={fill=white}, mark=*] table [x index=0, y index=1, skip coords between index={0}{10}]{\figfolderhomo KD_0.out};
\label{fig_homo_KD_0}
\addplot[very thick, color=colorB, mark options={fill=white}, mark=square*] table [x index=0, y index=3, skip coords between index={0}{10}]{\figfolderhomo KD_0.out};
\label{fig_homo_KD_1}
\addplot[very thick, color=colorC, mark options={fill=white}, mark=triangle*] table [x index=0, y index=5, skip coords between index={0}{10}]{\figfolderhomo KD_0.out};
\label{fig_homo_KD_8}
\addplot[very thick, dashed, black, domain=10:200,samples=2] {.58/x^(3/2)};
\label{fig_homo_KD_high_bound}
\node [draw=none,fill=none,anchor=south] at (axis cs: .05,0.001) {\shortstack[l]{
\ref{fig_homo_KD_high_bound} Eq. (\zelfhomofigboundB)\\
\vphantom{1}\\
\ref{fig_homo_KD_0} $\frac{\VAscaleparam\VAbackgroundspeedmeanflow-\VAtimestepvel_{{\VAtimestepno+1}}}{\sqrt{\VAbackgroundspeedtemp}}=0$ \\
\ref{fig_homo_KD_1} $\frac{\VAscaleparam\VAbackgroundspeedmeanflow-\VAtimestepvel_{{\VAtimestepno+1}}}{\sqrt{\VAbackgroundspeedtemp}}=1$ \\
\ref{fig_homo_KD_8} $\frac{\VAscaleparam\VAbackgroundspeedmeanflow-\VAtimestepvel_{{\VAtimestepno+1}}}{\sqrt{\VAbackgroundspeedtemp}}=8$}};
\end{axis}
\end{tikzpicture}%
}
\end{minipage}
  \captionof{figure}{The Wasserstein distance between the distribution of the kinetic process and of the KD process as a function of collisionality for a single time step. The experiment was conducted by averaging over 4,000,000 simulated particles and with $\VAtimestep=1$.}
  \label{fig:ap_liminf_W1err}
\end{figure}
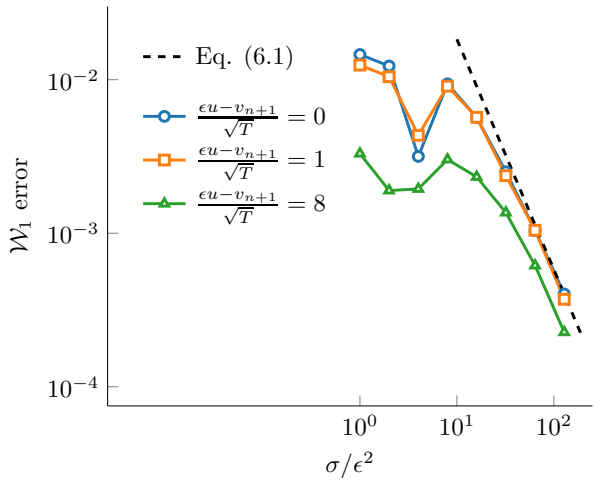

\paragraph{Numerical validation.} We illustrate the accuracy of the bound in Equation~\eqref{eq:ap_liminf_KD_distrerror} in Figure~\ref{fig:ap_liminf_W1err} with $\VAbackgroundspeedtemp=1$ and $\VAbackgroundspeedmeanflow=0$, conditioned on different dedimensionalized values for the velocity $\VAtimestepvel_{\VAtimestepno+1}$ at time $t^{(\VAtimestepno+1)}$. Independently of the value of the final velocity, the bound of Equation~\eqref{eq:ap_liminf_KD_distrerror} holds.

\subsection{Error propagation and total error in the diffusive limit\label{subsec:ap_liminf_KD_propagation}}

The second term on the right hand side of Equation~\eqref{eq:ap_liminf_errorpartsW1} has been treated in Section~\ref{subsec:ap_liminf_KD_distrerror}. In this Section, we will first treat the first term, expressing how the error at the start of the time step propagates to the next time step. Then, we will solve the resulting recursion to find the total error as in Section~\ref{subsec:ap_limnul_KD_propagation}.

An important distinction with the error analysis for fixed values of $\VAscaleparam$ is that now the kinetic contributions have a negligible effect. This means the conditioning on the velocities is not necessary here, and we can immediately use that propagation with the same operator can not increase the error. The first term of Equation~\eqref{eq:ap_liminf_errorpartsW1} can thus be bounded as
\begin{equation}
\zelfkalkdun{f(t^{(\VAtimestepno)})\!+\!\VAapKoperatoronlyD_{\VAtimestep}(f(t^{(\VAtimestepno)}))\!-\!\left(f^{(\VAtimestepno)}
\!+\!\VAapKoperatoronlyD_{\VAtimestep}(f^{(\VAtimestepno)})\right)}\leq\zelfkalkdun{f(t^{(\VAtimestepno)})-f^{(\VAtimestepno)}}\,.\label{eq:ap_liminf_KD_errprop}
\end{equation}

The bound in Equation~\eqref{eq:ap_liminf_KD_errprop} for the error propagation term in Equation~\eqref{eq:ap_liminf_errorpartsW1} together with the bound from Equation~\eqref{eq:ap_liminf_KD_distrerror} for the local error term in Equation~\eqref{eq:ap_liminf_errorpartsW1}, gives us the recursion formula
\begin{equation}
\zelfkalkdun{f(t^{(\VAtimestepno+1)})-f^{(\VAtimestepno+1)}}\leq\zelfkalkdun{f(t^{(\VAtimestepno)})-f^{(\VAtimestepno)}}+0.58\sqrt{\VAbackgroundspeedtemp}\frac{\VAscaleparam^{3}}{\sqrt{\VAratespatial^{3}\VAtimestep}}\,,\quad\VAscaleparam\rightarrow0\,.
\end{equation}
The solution of this recursion at the final simulation time $\VAtotaltime=\VAnoftimesteps\VAtimestep=t^{(\VAnoftimesteps)}$ is a bound on the Wasserstein error at that time, and equals
\begin{equation}
\zelfkalkdun{f(t^{(\VAnoftimesteps)})-f^{(\VAnoftimesteps)}}\leq0.58\sqrt{\VAbackgroundspeedtemp}\frac{\VAscaleparam^3}{\sqrt{\VAratespatial^3\VAtimestep^3}}\VAtotaltime\,,\qquad\VAscaleparam\rightarrow0\,.\label{eq:ap_liminf_KD_total}
\end{equation}
The error bound in Equation~\eqref{eq:ap_liminf_KD_total} proves the asymptotic preserving property, since, as the scale parameter becomes zero, the error due to replacing the real simulation with the KD scheme vanishes.

\section{Numerical experiment\label{sec:ap_num}}

In Sections~\ref{sec:ap_limnul} and~\ref{sec:ap_liminf}, it was proven that the simulation error by using the KD simulation of Section~\ref{sec:ap_newscheme}, is low in both the low-collisional and high-collisional regime. In this section, we numerically illustrate and extend these results with the simulation outcome of a model problem. We furthermore numerically show the, in Section~\ref{sec:ap_newscheme} predicted, computational speed-up in high-collisional regimes.

The model problem under study, is the Boltzmann equation of Equation~\eqref{eq:ap_kinetic_integrodiff} for spatially homogeneous parameters $\VAbackgroundspeedmeanflow(x)\equiv0$, $\VAbackgroundspeedtemp(x)\equiv1$, and $\VAratespatial(x)\equiv1$, for different values of $\VAscaleparam$. The initial condition is
\begin{equation}
\VAnsource({x},{v})=\delta_0(x)\left(\frac{1}{2\sqrt{2\pi}}e^{-\frac{(v+10)^2}{2}}+\frac{1}{2\sqrt{2\pi}}e^{-\frac{(v-10)^2}{2}}\right)\,,
\end{equation}
i.e., the particles start at $x=0$ and half of them have a normally distributed velocity with mean $-10$ and variance 1, and the other half have a normally distributed velocity with mean $10$ and variance $1$. The considered result of this model problem is the distribution at time $\VAtotaltime=1$ and for the KD simulation a single timestep of duration $\VAtimestep=\VAtotaltime=1$ is used.

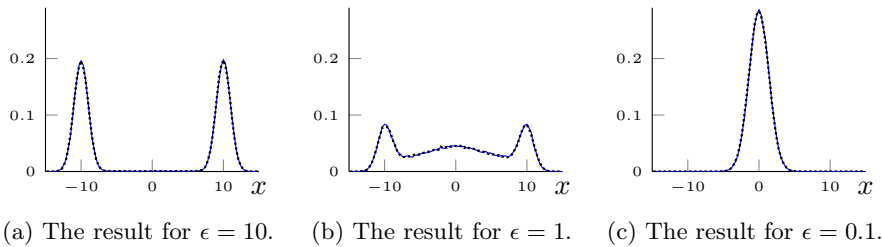
\begin{figure}[H]
\centering
\begin{subfigure}{.3\textwidth}
  \centering
  \resizebox{\textwidth}{!}{\begin{tikzpicture}
\begin{axis}[
  width=120pt,
  height=103pt,
  xmin=-15, xmax=15,
  ymin=0, ymax=.29,
  axis y line*=left,
  ytick={0,.1,.2},
  axis x line*=bottom,
  xlabel=$x$,
  xtick={-10,0,10},
  tick label style={font=\tiny},
  every axis x label/.style={at={(current axis.right of origin)},anchor=north},
  legend cell align={left},
  legend style={draw=none, fill=none,at={(axis cs:10,0.000002)}}]
\addplot[blue] table[x index=0, y index=1] {figuren/histogram_K_E10_P100000.txt};
\label{fig_res_K}
\addplot[thick, densely dotted, black] table[x index=0, y index=1] {figuren/histogram_KD_E10_P100000.txt};
\label{fig_res_KD}
\end{axis}
\end{tikzpicture}}
\caption{The result for $\VAscaleparam=10$.}
\label{fig:ap_num_hist_eps10}
\end{subfigure}\quad
\begin{subfigure}{.3\textwidth}
  \centering
  \resizebox{\textwidth}{!}{\begin{tikzpicture}
\begin{axis}[
  width=120pt,
  height=103pt,
  xmin=-15, xmax=15,
  ymin=0, ymax=.29,
  axis y line*=left,
  ytick={0,.1,.2},
  axis x line*=bottom,
  xlabel=$x$,
  xtick={-10,0,10},
  tick label style={font=\tiny},
  every axis x label/.style={at={(current axis.right of origin)},anchor=north},
  legend cell align={left},
  legend style={draw=none, fill=none,at={(axis cs:10,0.000002)}}]
\addplot[blue] table[x index=0, y index=1] {figuren/histogram_K_E1_P100000.txt};
\label{fig_res_K}
\addplot[densely dotted, black, thick]  table[x index=0, y index=1] {figuren/histogram_KD_E1_P100000.txt};
\label{fig_res_KD}
\end{axis}
\end{tikzpicture}}
\caption{The result for $\VAscaleparam=1$.}
\label{fig:ap_num_hist_eps1}
\end{subfigure}\quad
\begin{subfigure}{.3\textwidth}
  \centering
  \resizebox{\textwidth}{!}{\begin{tikzpicture}
\begin{axis}[
width=120pt,
  height=103pt,
  xmin=-15, xmax=15,
  ymin=0, ymax=.29,
  axis y line*=left,
  ytick={0,.1,.2},
  axis x line*=bottom,
  xlabel=$x$,
  xtick={-10,0,10},
  tick label style={font=\tiny},
  every axis x label/.style={at={(current axis.right of origin)},anchor=north},
  legend cell align={left},
  legend style={draw=none, fill=none,at={(axis cs:10,0.000002)}}]
\addplot[blue] table[x index=0, y index=1] {figuren/histogram_K_Ek1_P100000.txt};
\label{fig_res_K}
\addplot[densely dotted, black, thick]  table[x index=0, y index=1] {figuren/histogram_KD_Ek1_P100000.txt};
\label{fig_res_KD}
\end{axis}
\end{tikzpicture}}
\caption{The result for $\VAscaleparam=0.1$.}
\label{fig:ap_num_hist_epsk1}
\end{subfigure}
\caption{The probability distribution of the final position of the KD simulation (\ref{fig_res_KD}) and a standard MC method (\ref{fig_res_K}) for a model problem for different values of $\VAscaleparam$. For both methods $100,\!000$ particles were used.}
\end{figure}

The distribution of the position at time $\VAtotaltime$ resulting from the KD simulation is compared with the result of the standard MC method in Figures~\ref{fig:ap_num_hist_eps10}--\ref{fig:ap_num_hist_epsk1} for this model problem with $\VAscaleparam\in\{10,1,0.1\}$. In all three cases, the match is near-perfect. The values $\VAscaleparam=10$ and $\VAscaleparam=0.1$ correspond to respectively a low-collisional and high-collisional regime, and the close match is thus in accordance with the analytical results of Section~\ref{sec:ap_limnul}, respectively Section~\ref{sec:ap_liminf}. For the intermediate value $\VAscaleparam=1$, the close match in Figure~\ref{fig:ap_num_hist_eps1}, illustrates that the low error persists, even beyond the asymptotic cases studied in the previous two sections.

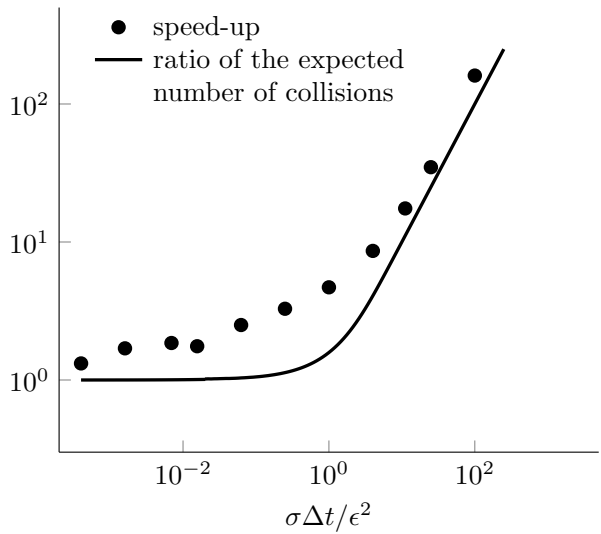
\begin{figure}[H]
\centering
\begin{minipage}{.65\textwidth}
  \centering
  \resizebox{\textwidth}{!}{\begin{tikzpicture}
\begin{axis}[
  width=240pt,
  height=207pt,
  xmin=2e-4, xmax=5e3,
  ymin=3e-1, ymax=5e2,
  axis y line*=left,
  ytick = {.0001, .01,1,10,100,10000},
  xtick = {.0001,.01,1,100,10000},
  axis x line*=bottom,
  xmode=log,
  ymode=log,
  xlabel=$\VAratespatial\VAtimestep/\VAscaleparam^2$,
  xlabel near ticks,
  legend cell align={left},
  legend style={draw=none, fill=none,at={(axis cs:10,0.000002)}}]
\addplot[draw=white, very thick, mark=*, mark options={black}] table[x index=0, y index=3] {figuren/speedup_P50000.txt};
\label{fig_speedup}
\addplot[very thick, black, domain=.0004:.0202,samples=100] {1/(1-x/2+x*x/6)};
\addplot[very thick, black, domain=.02:250,samples=100] {x/(1-exp(-x))};
\label{fig_speedup_theoretisch}
\node [draw=none,fill=none,anchor=south] at (axis cs: .08,90) {\shortstack[l]{
\ref{fig_speedup} speed-up\\
\ref{fig_speedup_theoretisch} ratio of the expected\\
\hphantom{\ref{fig_speedup_theoretisch}} number of collisions}};
\end{axis}
\end{tikzpicture}}
\end{minipage}
\caption{The speed-up by using the KD scheme instead of a standard MC method as function of the collisionality $\frac{\VAratespatial\VAtimestep}{\VAscaleparam^2}$. The results are found from numerical experiments with $50,\!000$ particles.}
\label{fig:ap_speedup}
\end{figure}

The motivation to use the KD scheme lies in its huge reduction of computational cost compared to the standard MC method when the collisionality is high. Then, most of the collisions that occur in the standard MC method, are not explicitly executed in the KD scheme, but aggregated into a diffusive step. The resulting speed-up in computational time is numerically computed and shown in Figure~\ref{fig:ap_speedup}. The computational speed-up is approximately proportional to the reduction in the number of collisions that has to be executed, as can also be seen in Figure~\ref{fig:ap_speedup}. In a standard MC method, the expected number of collisions occurring in a time step $\VAtimestep$ equals $\frac{\VAratespatial\VAtimestep}{\VAscaleparam^2}$, whereas in a KD simulation only the first potential collision is actually executed. Hence, the expected number of collisions in a KD simulation equals the probability of at least one collision occurring in a standard MC simulation, which equals $1-e^{-\frac{\VAratespatial\VAtimestep}{\VAscaleparam^2}}$. The ratio of both is shown in Figure~\ref{fig:ap_speedup}.

\section{Conclusion}

To overcome high simulation costs in high-collisional regimes and erroneous results in low-collisional regimes, we present a hybridized simulation algorithm for the Boltzmann-BGK equation that combines the standard kinetic simulation and a random walk simulation that captures the high-collisional behaviour, at a much reduced computational cost. We have proven consistency with the kinetic process and we have shown the asymptotic behaviour in the high-collisional limit is also preserved. Furthermore, resulting from the matching first two moments, errors are expected to be low for intermediate collisionality, as is numerically illustrated.

To keep the exposition simple, we have ignored absorption in this paper. Absorption can be easily included in the scheme by simulating the time to the next absorption separately from the other motion. If an absorption would then occur during a diffusive substep, the full diffusive substep can be replaced with kinetic motion or it can be replaced by a diffusive substep until the absorption time.

An important improvement has been presented in~\cite{mortier2019KDfusioncase}, where the scheme was applied to a realistic fusion related case. There, an advection term was added to cope with the heterogeneity of the collision rate, and the treatment of reflective boundary conditions was presented.

In future work, we will focus further on the application of this simulation scheme to fusion reactor simulations, where neutrals play a crucial role in the plasma edge by providing source terms for the plasma equations. The neutrals undergo both low-collisional and high-collisional regimes, motivating the use of hybridized schemes. The main challenge in applying the algorithm presented here in that context, is the extraction of source terms from the diffusion part of the simulation. Other possible extensions include the addition of a deterministic fluid model as control variate and multilevel Monte Carlo where the time step size can be used as a level parameter.

\section*{Acknowledgments}

The first author was funded a personal grant from the Research Foundation - Flanders (FWO) under fellowship number 1189919N.\\
This work has been carried out within the framework of the EUROfusion Consortium and has received funding from the Euratom research and training programme 2014-2018 and 2019-2020 under grant agreement No 633053. The views and opinions expressed herein do not necessarily reflect those of the European Commission.


\bibliographystyle{abbrv}

\bibliography{referenties_doctoraat_kopie}

\end{document}